\documentclass[journal,twoside,letterpaper]{IEEEtran}

\usepackage{amssymb,amsthm,amsmath}
\usepackage[colorlinks=true]{hyperref}
\usepackage{nicefrac,enumerate,cite}
\usepackage[capitalise,sort]{cleveref}
\usepackage{graphicx}

\numberwithin{equation}{section}
\newtheorem{lemma}{Lemma}
\newtheorem{remark}[lemma]{Remark}
\newtheorem{corollary}[lemma]{Corollary}
\newtheorem{definition}[lemma]{Definition}
\newtheorem{theorem}[lemma]{Theorem}
\newtheorem{proposition}[lemma]{Proposition}
\newtheorem{example}[lemma]{Example}

\newcommand{\R}{\mathbb{R}}
\newcommand{\N}{\mathbb{N}}

\renewcommand{\epsilon}{\varepsilon}

\newcommand{\ssum}[2]{\mathop{\textstyle{\sum}}_{#1}^{#2}}
\newcommand{\sprod}[2]{\mathop{\textstyle{\prod}}_{#1}^{#2}}

\newcommand{\mednorm}[1]{\| #1 \|}
\newcommand{\opnorm}[1]{\mednorm{#1}}

\newcommand{\ID}[1]{\mathrm{id}_{\R^{#1}}}

\newcommand{\auxproofcomplayerreal}[1]{G_{#1}}

\newcommand{\sqfct}{\mathrm{sq}}
\newcommand{\prodfct}{\mathrm{pr}}
\newcommand{\maxfct}[1]{\mathrm{max}_{#1}}
\newcommand{\expfct}{\mathrm{e}}

\newcommand{\ANNs}{\mathcal{N}}
\newcommand{\neuronsANN}[1]{\mathcal{P}(#1)}
\newcommand{\depthANN}[1]{\mathcal{D}(#1)}

\newcommand{\DimANN}[2]{l_{#1}^{#2}}
\newcommand{\inDimANN}[1]{\DimANN{0}{#1}}
\newcommand{\outDimANN}[1]{\DimANN{\depthANN{#1}}{#1}}

\newcommand{\concANN}{\circ}

\newcommand{\affANN}[2]{\mathcal{A}_{#2}^{#1}}
\newcommand{\fctANN}[1]{\mathcal{R}_a^{#1}}

\newcommand{\IDANN}[1]{I_{#1}}
\newcommand{\idconstANN}{c}

\newcommand{\paralANN}[2]{\paralnbANN{#1}(#2)}
\newcommand{\paralnbANN}[1]{p_{#1}}

\newcommand{\ReLU}{\mathrm{ReLU}}
\newcommand{\fctReLUANN}[1]{\mathcal{R}_{\ReLU}^{#1}}

\newcommand{\Catalog}{\mathcal{F}}

\newcommand{\CatDimD}[1]{\mathcal{I}(#1)}
\newcommand{\CatDimT}[1]{\mathcal{O}(#1)}

\newcommand{\Lip}{L}

\newcommand{\CostLipwo}[2]{\mathrm{Cost}_{#1,#2}}

\newcommand{\CostApprox}[6]{\CostLipwo{#1}{#2}(#3,#4,#5,#6)}

\newcommand{\ApproxCat}[2]{[a,\weight,\approxSet{},\Lip,#1,#2]}

\newcommand{\ApproxExCat}[3]{[\ReLU,#3,\approxSet{},\Lip,#1,#2]}

\newcommand{\Lipconst}{K}
\newcommand{\Lipcatalog}{\mathcal{F}_{\Lipconst}^{\mathrm{Lip}}}
\newcommand{\Lipmaxcatalog}{\mathcal{F}_{\Lipconst}^{\mathrm{Lip,max}}}
\newcommand{\Lipprodcatalog}{\mathcal{F}_{\Lipconst}^{\mathrm{Lip,prod}}}
\newcommand{\Lipprod}{\mathcal{F}^{\mathrm{prod}}}
\newcommand{\GaussianRBFcat}{\mathcal{F}^{\mathrm{RBF}}}

\newcommand{\CNfixed}{\mathcal{C}_{\Catalog}^{l_0,\dots,l_{2D}}}
\newcommand{\CNarb}{\mathcal{C}_{\Catalog}}

\newcommand{\CNdepth}[1]{\mathcal{D}_{#1}}
\newcommand{\CNinput}[1]{\mathcal{I}_{#1}}
\newcommand{\CNoutput}[1]{\mathcal{O}_{#1}}
\newcommand{\CNwidth}[1]{\mathcal{W}_{#1}}
\newcommand{\CNtransl}[1]{\mathcal{B}_{#1}}

\newcommand{\CNelt}{\xi}
\newcommand{\CNelementfixed}{[(V_1,b_1, \allowbreak (f_{1,1},\dots,f_{1,n_1})),  \allowbreak \dots, \allowbreak (V_D,b_D, \allowbreak (f_{D,1},\dots,f_{D,n_D}))]}

\newcommand{\CNaff}[2]{\mathcal{A}^{#2,#1}}
\newcommand{\CNnonlin}[2]{\mathcal{G}^{#2,#1}}
\newcommand{\CNfct}[1]{\mathcal{R}^{#1}}

\newcommand{\CNrdlayer}[2]{Q^{#2,#1}}
\newcommand{\CNredom}[1]{\mathrm{Dom}_{\approxSet{} , #1}}

\newcommand{\CNlayerLip}[2]{\Lip^{#2,#1}}
\newcommand{\CNfullLip}[1]{\mathrm{Lip}_{\Lip , #1}}
\newcommand{\CNtechLip}[1]{\mathcal{T}_{\Lip , #1}}

\newcommand{\weight}{w}
\newcommand{\approxSet}[1]{Q_{#1}}

\newcommand{\OoGone}{{s_1}}
\newcommand{\OoGtwo}{{s_2}}

\title{\LARGE{Efficient Approximation of High-dimensional \\ Functions With Neural Networks}}

\author{Patrick Cheridito \quad Arnulf Jentzen \quad Florian Rossmannek\thanks{P. Cheridito and F. Rossmannek are with the Department of Mathematics at ETH Zurich. A. Jentzen is with the Faculty of Mathematics and Computer Science at the University of M\"unster.}}

\date{2021}

\begin{document}

\maketitle

\begin{abstract}
In this paper, we develop a framework for showing that neural networks 
can overcome the curse of dimensionality in different high-dimensional approximation problems.
Our approach is based on the notion of a catalog network, which is a generalization of a
standard neural network in which the nonlinear activation functions can vary from layer to layer as long 
as they are chosen from a predefined catalog of functions. As such, catalog networks 
constitute a rich family of continuous functions. We show that under appropriate conditions on the catalog, 
catalog networks can efficiently be approximated with rectified linear unit-type networks and provide precise estimates on 
the number of parameters needed for a given approximation accuracy. As special cases of the general results,
we obtain different classes of functions that can be approximated with ReLU networks without the curse 
of dimensionality.

\vskip 2mm

\end{abstract}


\section{Introduction}
\label{sec_intro}

\IEEEPARstart{M}{any} classical numerical approximation schemes that work well in low dimensions suffer from the so 
called curse of dimensionality, meaning that to achieve a desired approximation accuracy, their 
complexity has to grow exponentially in the dimension. On the other hand, neural networks have shown 
remarkable performance in different high-dimensional approximation problems. 
In this paper we prove that different classes of high-dimensional functions admit a 
neural network approximation without the curse of dimensionality. To do that, we introduce the notion of a catalog network, 
which is a generalization of a standard feedforward neural network in which the nonlinear activation functions can 
vary from one layer to another as long as they are chosen from a given catalog of continuous functions. 
We first study the approximability of different catalogs with neural networks. Then we show how the approximability 
of a catalog translates into the approximability of the corresponding catalog networks.
An important building block of our proofs is a new way of parallelizing networks that saves parameters compared 
to the standard parallelization. As special cases of our general results we obtain that different combinations of 
one-dimensional Lipschitz functions, sums, maxima and products as well as certain ridge functions and 
generalized Gaussian radial basis function networks admit a neural network approximation
without the curse of dimensionality. 

It has been shown that neural networks with a single hidden layer can approximate any finite-dimensional 
continuous function uniformly on compact sets arbitrarily well if they are allowed to have sufficiently many hidden neurons;
see, e.g., \cite{Funahashi1989,Cybenko1989, HorStinWhite1989, Hornik1991, LeshnoLinPinkusScho1993, GuliIsm2018b}.
Moreover, \cite{Jones1992} and \cite{Barron1993} have proved an $\mathcal{O}(n^{-1/2})$-rate for approximating 
functions in the $L^2$-norm with single-hidden-layer sigmoidal networks with $n$ neurons.
In particular, this breaks the curse of dimensionality, but it only applies to a special class of functions.
Since then, their results have been applied and generalized in different directions, always yielding rates of the same nature, 
but always applicable only to similarly restricted classes of functions.
For example, in \cite{DonGurvDarSont1997,GirAnz1993} these results have been extended to the $L^p$-norm for $1 \leq p < \infty$ and $p=\infty$, respectively, and
in \cite{KurSang2008} the approximation rate has been improved to a geometric rate for single functions.
However, the basis $\delta$ of the geometric rate $\mathcal{O}((1-\delta)^n)$ is usually not known. So, it could be so small
that the geometric rate does not give useful bounds for typical sizes of $n$. For further generalizations, see, e.g.,
\cite{Barron1992,Barron1994,GurvKoir1997,KlusBar2018,KurKaiKre1997,
KurSang2002,Kurkova2008,KaiKurSang2009,KaiKurSang2012}.
All of them use single-hidden-layer networks. However, neural networks with more than one hidden layer have
shown better performance in a number of applications; see, e.g., \cite{LeCun2015, GoodfBengioCourv2016} and the 
references therein. This has also been supported by theoretical evidence; for instance, in \cite{EldanShamir2016},
an example of a simple continuous function on $\R^d$ has been given that is expressible as a small feedforward network 
with two hidden layers but cannot be approximated with a single-hidden-layer network 
to a given constant accuracy unless its width is exponential in the dimension. Similarly,
it has been shown in \cite{SafranShamir2017} that indicator functions of $d$-dimensional balls 
can be approximated much more efficiently with two hidden layers than with one.
Related results for functions on the product of two $d$-dimensional spheres have been 
provided by \cite{Daniely2017}. 

\cite{Pinkus1999,MaioPinkus1999,GuliIsm2018a} have constructed special activation functions
which, in principle, allow to approximate every continuous function $f \colon [0,1]^d \to \R$ 
to any desired precision when used in a two-hidden-layers network with as few as $d$ neurons in the first and 
$2d + 2$ neurons in the second hidden layer. Theoretically, this breaks the curse of dimensionality quite 
spectacularly. However, it can be shown that the approximation result only holds if the size of the network 
weights is allowed to grow faster than polynomially in the inverse of the approximation error; see, e.g., 
\cite{BolcGrohsKutynPete2019,PeteVoigt2018}.

Further studies of the approximation capacity of neural networks with standard activation functions include, e.g,
\cite{Mhaskar1993,VoigtPete2019,PeteVoigt2018,LuShenYangZhang2020,Yarotsky2017}.
Their approach is based on approximating functions with polynomials and then approximating these polynomials with neural networks.
Polynomials can approximate smooth functions reasonably well, and neural networks are known to be able to approximate 
monomials efficiently. However, since the number of monomials needed to generate all polynomials in $d$ variables of order 
$k$ is $\binom{k+d}{d}$, the intermediate step from monomials to polynomials introduces a curse of dimensionality. It has been 
shown in \cite{Yarotsky2017} that this cannot be side-stepped. For instance, it is provably impossible 
to approximate the unit ball in the Sobolev space of any regularity with ReLU networks without the curse of dimensionality.
So, to overcome the curse of dimensionality with ReLU networks, one has to concentrate on special classes of functions.
\cite{Jones1992,Barron1993} and their extensions offer one such class. 
Coming from a different angle, \cite{MhasMicc1994} has obtained the same rate for
periodic functions with an absolutely convergent Fourier series.
In \cite{SchwabZech2019} the approximability of ``separately holomorphic" maps via Taylor expansions
and applications to parametric PDEs have been studied. The approach of \cite{SchwabZech2019} is again 
based on the intermediate approximation of polynomials, but the holomorphy ensures that the 
approximating polynomials contain only few monomials. \cite{GrohsHorJenWurs2018,JenSaliWelti2018,HutzJenKruseNguyen2020}
have proved that solutions of various PDEs admit neural network approximations without the curse of dimensionality.
Their arguments use the hierarchic structure of neural networks, which has more extensively been exploited in 
\cite{LeeGeMaRistArora2017,BolcGrohsKutynPete2019,GrohsPerekElbBolc2019}.
These papers are similar in spirit to ours since they also start from a ``basis'' of functions, which they 
approximate with neural networks and then use to build more complex functions.
However, \cite{BolcGrohsKutynPete2019,GrohsPerekElbBolc2019} do not study approximation rates in terms of the dimension.
On the other hand, in \cite{LeeGeMaRistArora2017} the curse of dimensionality is overcome, but the ``basis" in 
\cite{LeeGeMaRistArora2017} consists of the functions considered in \cite{Barron1993}. 
In this paper we consider more explicit classes of functions and 
provide bounds on the number of parameters needed
to approximate $d$-dimensional functions up to accuracy $\varepsilon$.

The rest of the paper is organized as follows. In Section \ref{section_DNNs}, we first establish the notation.
Then we recall basic facts from \cite{PeteVoigt2018,GrohsHorJenZimm2019,JenSaliWelti2018} on concatenating and 
parallelizing neural networks before we introduce a new way of network parallelization. 
In Section \ref{section_CNs}, we introduce the concepts of an approximable catalog and 
a catalog network. Section \ref{section_ex_approx_cat} is devoted to different
concrete examples of catalogs and a careful study of their approximability. In
Sections \ref{section_approx_results} and \ref{section_log_modification}, we derive bounds on the 
number of parameters needed to approximate a given catalog network to a desired accuracy with neural networks.
\cref{theorem_main,theorem_main_LOG} are the main results of this article.
In Section \ref{section_applications}, we derive different classes of high-dimensional 
functions that are approximable with ReLU networks without the curse of dimensionality.
Section \ref{section_conclusion} concludes. All proofs are relegated to the Appendix.


\section{Notation and Preliminary Results}
\label{section_DNNs}

A neural network encodes a succession of affine and non-linear transformations. Let us
denote $\N = \{1,2,\dots\}$ and consider the set of neural network skeletons
\begin{equation*}
	\ANNs = \bigcup_{D \in \N} \bigcup_{(l_0,\dots,l_D) \in \N^{D+1}} \prod_{k=1}^{D} ( \R^{l_k \times l_{k-1}} \times \R^{l_k} ).
\end{equation*}
We denote the depth of a neural network skeleton $\phi \in \ANNs$ by $\depthANN{\phi} = D$, 
the number of neurons in the $k$th layer by $\DimANN{k}{\phi} = l_k$, $k \in \{0,\dots,D\}$,
and the number of network parameters by $\neuronsANN{\phi} = \sum_{k=1}^D l_k(l_{k-1}+1)$.
Moreover, if $\phi \in \ANNs$ is given by $\phi = [(V_1,b_1),\dots,(V_D,b_D)]$, we denote by 
$\affANN{\phi}{k} \in C(\R^{l_{k-1}},\R^{l_{k}})$, $k \in \{1,\dots,D\}$, the affine function $x \mapsto V_kx+b_k$.
Let $a \colon \R \to \R$ be a continuous activation function. As usual, we extend it, for every positive integer $d$,
to a function from $\R^d$ to $\R^d$ mapping $(x_1, \dots, x_d)$ to $(a(x_1), \dots, a(x_d))$.
Then the $a$-realization of $\phi \in \ANNs$ is the function $\fctANN{\phi} \in C(\R^{l_0},\R^{l_D})$ given by
\begin{equation*}
\fctANN{\phi} = \affANN{\phi}{D} \circ a \circ \affANN{\phi}{D-1} \circ \cdots a \circ \affANN{\phi}{1}.
\end{equation*}
We recall that suitable $\phi_1,\phi_2 \in \ANNs$ can be composed such that the $a$-realization of 
the resulting network equals the concatenation $\fctANN{\phi_2} \circ \fctANN{\phi_1}$.
This is done by combining the output layer of $\phi_1$ with the input layer of $\phi_2$.
More precisely, if $\phi_1 = [(V_1,b_1),\dots,(V_D,b_D)]$ and $\phi_2 = [(W_1,c_1),\dots,(W_E,c_E)]$ satisfy $\outDimANN{\phi_1} = \inDimANN{\phi_2}$, then the concatenation $\phi_2 \concANN \phi_1 \in \ANNs$ is given by
\begin{equation*}
\begin{split}
	\phi_2 \concANN \phi_1 = [&(V_1,b_1),\dots,(V_{D-1},b_{D-1}), \\
	&(W_1V_D,W_1b_D+c_1),(W_2,c_2),\dots,(W_E,c_E)].
\end{split}
\end{equation*}
The following result is straight-forward from the definition.
A formal proof can be found in \cite{GrohsHorJenZimm2019}.

\begin{proposition}
\label{prop_ANN_concat}
	The concatenation
\begin{equation*}
	(\cdot) \concANN (\cdot) \colon \{(\phi_1,\phi_2) \in \ANNs \times \ANNs \colon \outDimANN{\phi_1} = \inDimANN{\phi_2} \} \rightarrow \ANNs
\end{equation*}
is associative and for all $\phi_1,\phi_2 \in \ANNs$ with $\outDimANN{\phi_1} = \inDimANN{\phi_2}$ one has
\begin{enumerate}[\rm 1)]
\item $\fctANN{\phi_2 \concANN \phi_1} = \fctANN{\phi_2} \circ \fctANN{\phi_1}$ for all $a \in C(\R,\R)$,

\item $\depthANN{\phi_2 \concANN \phi_1} = \depthANN{\phi_1} + \depthANN{\phi_2} - 1$,

\item $\DimANN{k}{\phi_2 \concANN \phi_1} = \DimANN{k}{\phi_1}$ if $k \in \{0,\dots,\depthANN{\phi_1}-1\}$,

\item $\DimANN{k}{\phi_2 \concANN \phi_1} = \DimANN{k+1-\depthANN{\phi_1}}{\phi_2}$ if $k \in \{\depthANN{\phi_1},\dots,\depthANN{\phi_2 \concANN \phi_1}\}$,

\item\label{prop_ANN_concat_item_param} $\neuronsANN{\phi_2 \concANN \phi_1} = \neuronsANN{\phi_1} + \neuronsANN{\phi_2} + \DimANN{1}{\phi_2}\DimANN{\depthANN{\phi_1}-1}{\phi_1} - \inDimANN{\phi_2}\DimANN{1}{\phi_2} - \outDimANN{\phi_1}(\DimANN{\depthANN{\phi_1}-1}{\phi_1} + 1)$,

\item\label{prop_ANN_concat_item_depth_1_1} $\neuronsANN{\phi_2 \concANN \phi_1} \leq \neuronsANN{\phi_1}$ if $\depthANN{\phi_2} = 1$ and $\DimANN{1}{\phi_2} \leq \outDimANN{\phi_1}$

\item\label{prop_ANN_concat_item_depth_1_2} and $\neuronsANN{\phi_2 \concANN \phi_1} \leq \neuronsANN{\phi_2}$ if $\depthANN{\phi_1} = 1$ and $\inDimANN{\phi_1} \leq \inDimANN{\phi_2}$.
\end{enumerate}
\end{proposition}

The next lemma is a direct consequence of the above and will be used later to estimate the number of parameters in our approximating networks.

\begin{lemma}
\label{lemma_parameters_of_id_conc}
	Let $a \in C(\R,\R)$ and $\phi \in \ANNs$.
Suppose that $\psi_1,\psi_2 \in \ANNs$ satisfy $\depthANN{\psi_1} = \depthANN{\psi_2} = 2$, $\inDimANN{\psi_1} = \DimANN{2}{\psi_1} = \inDimANN{\phi}$ and $\inDimANN{\psi_2} = \DimANN{2}{\psi_2} = \outDimANN{\phi}$.
Denote $m = \DimANN{1}{\psi_1}$ if $\depthANN{\phi} = 1$ and $m = \DimANN{\depthANN{\phi}-1}{\phi}$ if $\depthANN{\phi} \geq 2$.
Then
\begin{equation*}
\begin{split}
	\neuronsANN{\psi_2 \concANN \phi \concANN \psi_1} &= \neuronsANN{\phi} + \DimANN{1}{\psi_1} (\inDimANN{\phi}+1) + \DimANN{1}{\psi_2} (\DimANN{\depthANN{\phi}}{\phi}+1) \\
	&\quad + \DimANN{1}{\phi} (\DimANN{1}{\psi_1} - \inDimANN{\phi}) + m ( \DimANN{1}{\psi_2}  - \DimANN{\depthANN{\phi}}{\phi} ).
\end{split}
\end{equation*}
\end{lemma}

The standard parallelization of two network skeletons $\phi_1 = [(V_1,b_1),\dots,(V_D,b_D)]$ 
and $\phi_2 = [(W_1,c_1),\dots,(W_D,c_D)]$ of the same depth is given by $\paralANN{}{\phi_1,\phi_2} =$
\begin{equation*}
	\bigg[ \bigg(
	\begin{bmatrix}
		V_1 & 0 \\
		0 & W_1
	\end{bmatrix},
	\begin{bmatrix}
		b_1 \\
		c_1
	\end{bmatrix} \bigg), \dots, \bigg(
	\begin{bmatrix}
		V_D & 0 \\
		0 & W_D
	\end{bmatrix},
	\begin{bmatrix}
		b_D \\
		c_D
	\end{bmatrix} \bigg) \bigg].
\end{equation*}
From there, arbitrarily many network skeletons $\phi_1,\dots,\phi_n \in \ANNs$, $n \in \N_{\geq 3}$, of the same depth 
can be parallelized iteratively:
\begin{equation*}
	\paralANN{}{\phi_1,\dots,\phi_n} = \paralANN{}{\paralANN{}{\phi_1,\dots,\phi_{n-1}},\phi_n}.
\end{equation*}
The first three statements of the next proposition follow immediately from the definition.
The last one is shown in \cite{GrohsHorJenZimm2019}.

\begin{proposition}
\label{prop_ANN_paral_same_depth}
	The parallelization
\begin{equation*}
	\paralnbANN{} \colon \bigcup_{n \in \N} \{ (\phi_1,\dots,\phi_n) \in \ANNs^n \colon \depthANN{\phi_1} = \dots = \depthANN{\phi_n} \} \rightarrow \ANNs
\end{equation*}
satisfies for all $\phi_1,\dots,\phi_n \in \ANNs$, $n \in \N$, with the same depth
\begin{enumerate}[\rm 1)]

\item $\fctANN{\paralANN{}{\phi_1,\dots,\phi_n}}(x_1,\dots,x_n) = (\fctANN{\phi_1}(x_1),\dots,\fctANN{\phi_n}(x_n))$ for all $x_1 \in \R^{\inDimANN{\phi_1}}, \dots, x_n \in \R^{\inDimANN{\phi_n}}$ and each $a \in C(\R,\R)$,

\item\label{prop_ANN_paral_same_depth_item} $\DimANN{k}{\paralANN{}{\phi_1,\dots,\phi_n}} = \sum_{j=1}^{n} \DimANN{k}{\phi_j}$ for all $k \in \{0,\dots,\depthANN{\phi_1}\}$,

\item\label{prop_ANN_paral_same_depth_item_same_network} $\neuronsANN{\paralANN{}{\phi_1,\dots,\phi_n}} \leq n^2 \neuronsANN{\phi_1}$ whenever $\DimANN{k}{\phi_i} = \DimANN{k}{\phi_j}$ for all $k \in \{0,\dots,\depthANN{\phi_1}\}$ and all $i,j \in \{1,\dots,n\}$

\item\label{prop_ANN_paral_same_depth_item_quadratic_param} and $\neuronsANN{\paralANN{}{\phi_1,\dots,\phi_n}} \leq \frac{1}{2} \big[ \sum_{j=1}^{n} \neuronsANN{\phi_j} \big]^2$.

\end{enumerate}
\end{proposition}

Neural networks with different depths can still be parallelized, but only for a special class of activation functions.

\begin{definition}
\label{def_c_id_requ}
We say a function $a \in C(\R,\R)$ fulfills the $\idconstANN$-identity requirement for a number 
$\idconstANN \ge 2$ if there exists $\IDANN{} \in \ANNs$ such that $\depthANN{\IDANN{}} = 2$, $\DimANN{1}{\IDANN{}} \leq \idconstANN$ and $\fctANN{\IDANN{}} = \ID{}$.
\end{definition}

Note that if $I$ satisfies $\fctANN{\IDANN{}} = \ID{}$, one can also realize the identity function 
$\ID{d}$ for any $d \in \N$, using $d$-fold parallelization $\IDANN{d} = \paralANN{}{\IDANN{},\dots,\IDANN{}}$.
Obviously, $\DimANN{1}{\IDANN{d}} \leq \idconstANN d$.

The most prominent example satisfying \cref{def_c_id_requ} is the rectified linear unit activation 
$\ReLU \colon \R \rightarrow \R, ~x \mapsto \max\{x,0\}$.
It fulfills the $2$-identity requirement with $\IDANN{} = [ ([1 \ -1]^T,[0 \ 0]^T),$ $([1 \ -1],0) ]$.
However, it is easy to see that generalized ReLU functions of the form 
\begin{equation*}
	a(x) =
	\begin{cases}
	rx & \mbox{ if } x \geq 0\\
	sx & \mbox{ if } x < 0
	\end{cases}
\end{equation*}
for $(r,s) \in \R^2$ with $r+s \ne 0$, such as leaky ReLU, also satisfy the $2$-identity requirement.\footnote{Other 
activation functions satisfying the identity requirement are polynomials. For example, $\frac{1}{2}((x+1)^2-x^2-1) = x$ 
shows this for $x^2$.}

\begin{figure}
\centering
\includegraphics{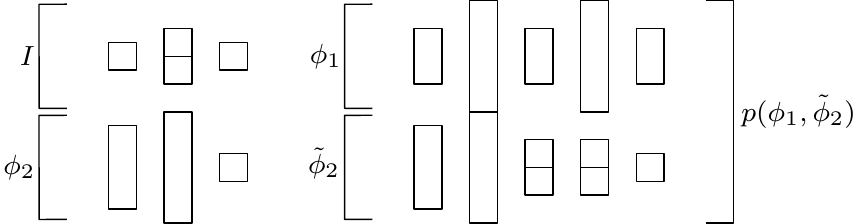}
\caption{Parallelization of a network $\phi_1$ (depth 4) and a shorter network $\phi_2$ (depth 2) obtained by concatenating $\phi_2$ twice with a network $I$ arising from the 2-identity requirement, resulting in $\tilde{\phi}_2 = \IDANN{} \concANN \IDANN{} \concANN \phi_2$.}
\label{fig_parallelization}
\end{figure}

Using the identity requirement, one can parallelize networks of arbitrary depths.
If $\phi_1,\dots,\phi_n \in \ANNs$ have different depths, one simply concatenates the 
shorter ones with identity networks until all have the same depth. Then one applies 
the standard parallelization; see Figure \ref{fig_parallelization} for an illustration.

\begin{figure}
\centering
\includegraphics{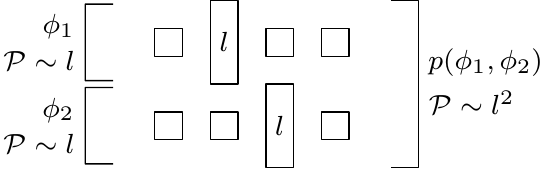}
\caption{The parallelization of $\phi_1$ with architecture $(1,l,1,1)$ and $\phi_2$ with architecture 
$(1,1,l,1)$ has more than $l^2$ parameters.}
\label{fig_parallelization_squaring}
\end{figure}

\begin{figure}
	\centering
	\includegraphics{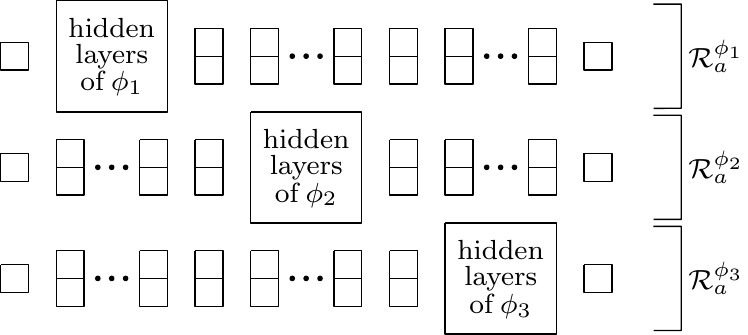}
	\caption{New ``diagonalized'' parallelization resulting from shifting $\phi_1,\phi_2,\phi_3$ away from each other.}
	\label{fig_parallelization_novel}
\end{figure}

Although this successfully parallelizes networks with arbitrary architecture, one can do better in terms of parameter counts.
The estimate in \cref{prop_ANN_paral_same_depth}.\eqref{prop_ANN_paral_same_depth_item_quadratic_param}
contains a square of $\sum_{j=1}^{n} \neuronsANN{\phi_j}$. This is not due to lax estimates, but a square can actually appear
if, for some $j$, there are two large consecutive layers in $\paralANN{}{\phi_j,\phi_{j+1}}$ which in $\phi_j$ and $\phi_{j+1}$
were next to small layers; see Figure \ref{fig_parallelization_squaring}. To avoid this, we introduce a new parallelization
which uses identity networks to shift $\phi_1, \dots, \phi_n$ away from each other and, as a result, achieves a 
parameter count that is linear in $\sum_{j=1}^{n} \neuronsANN{\phi_j}$. For instance, to parallelize $\phi_1$ and $\phi_2$,
we add $\depthANN{\phi_2}$ identity networks after $\phi_1$ and $\depthANN{\phi_1}$ identity networks in front of $\phi_2$ 
before applying $p$. The realization of the resulting network still is $(x_1,x_2) \mapsto (\fctANN{\phi_1}(x_1),\fctANN{\phi_2}(x_2))$.
Extending this construction to more than two networks is straight-forward; see Figure \ref{fig_parallelization_novel}.
We denote it by $\paralnbANN{\IDANN{}}$, where $\IDANN{} \in \ANNs$ is the network satisfying the identity requirement.
The following proposition shows that $\paralnbANN{\IDANN{}}$ achieves our goal of a linear parameter count in 
$\ssum{j=1}{n} \neuronsANN{\phi_j}$.

\begin{proposition}
\label{prop_ANN_paral_NEW}
Assume $a \in C(\R,\R)$ fulfills the $\idconstANN$-identity requirement for a number 
$\idconstANN \ge 2$ with $\IDANN{} \in \ANNs$.
Then the parallelization $\paralnbANN{\IDANN{}} \colon \bigcup_{n \in \N} \ANNs^n \rightarrow \ANNs$ satisfies 
\begin{equation*}
	\neuronsANN{\paralANN{\IDANN{}}{\phi_1,\dots,\phi_n}} \leq \big( \tfrac{11}{16} c^2 l^2 n^2 - 1 \big) \ssum{j=1}{n} \neuronsANN{\phi_j}
\end{equation*}
for all $n \in \N$ and $\phi_1,\dots,\phi_n \in \ANNs$, where we denote 
$l = \max_{j \in \{1,\dots,n\}} \max\{\inDimANN{\phi_j},\outDimANN{\phi_j}\}$.
\end{proposition}

It can be seen from the proof that the inequality of \cref{prop_ANN_paral_NEW} is never an equality.
However, it can be shown that it is asymptotically sharp up to a constant for large $n$.
Indeed, if $\idconstANN = 2$ (as is the case for $\ReLU$) and if $\phi_1 = \dots = \phi_n$ has depth at least two ($\depthANN{\phi_1} \geq 2$) and a single neuron in each layer ($\DimANN{k}{\phi_1} = 1$ for all $k$), then
\begin{equation*}
	(2n^3-n^2) \neuronsANN{\phi_1} = \neuronsANN{\paralANN{\IDANN{}}{\phi_1,\dots,\phi_n}} \leq \big( \tfrac{11}{4} n^3 - n \big) \neuronsANN{\phi_1}.
\end{equation*}
The first inequality is verified in the Appendix while the second one is a consequence of \cref{prop_ANN_paral_NEW}.
Hence, the bound in the proposition is asymptotically sharp up to a factor of at most $\frac{11}{8}$.

\begin{figure}
	\centering
	\includegraphics{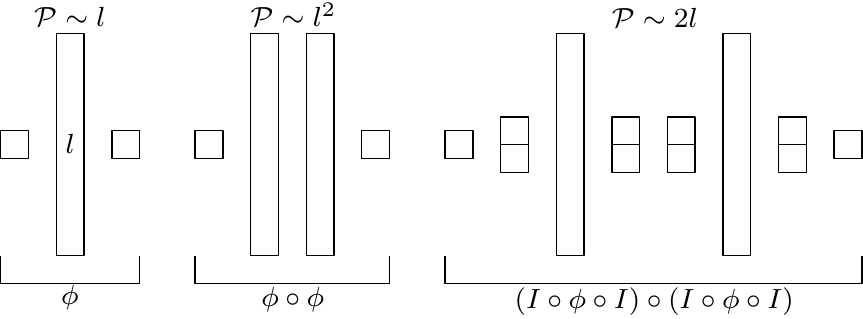}
	\caption{Concatenation with and without additional identity networks. Here, $\phi$ is a network of depth 2 
	with $l$ neurons in its hidden layer, and $I$ is assumed to satisfy the 2-identity requirement.}
	\label{fig_cheaper_concatenation}
\end{figure}

\cref{prop_ANN_paral_NEW} illustrates that there is a fundamental difference between counting the number of 
neurons and counting the number of parameters. As already observed in \cite{PeteVoigt2018,GrohsHorJenZimm2019,JenSaliWelti2018},
this also plays a role for the concatenation. The standard concatenation of two networks $\phi_1$ and $\phi_2$
has roughly ${\cal P}(\phi_1) + {\cal P}(\phi_2)$ neurons. But the parameter count 
may increase much more dramatically. If, e.g., most of the neurons of $\phi_1$ are in the last hidden layer and 
most of the neurons of $\phi_2$ in the first hidden layer,
then $\phi_2 \circ \phi_1$ has roughly $\neuronsANN{\phi_1} \cdot \neuronsANN{\phi_2}$ parameters; see Figure \ref{fig_cheaper_concatenation}. To counter this, one can use the concatenation 
\[
I_{l^{\phi_2}_{{\cal D}(\phi_2)}} \circ \phi_2 \circ I_{l^{\phi_2}_0}
\circ I_{l^{\phi_1}_{{\cal D}(\phi_1)}} \circ \phi_1 \circ I_{l^{\phi_1}_0}
\]
instead of $\phi_2 \circ \phi_1$, where $I_d$ is an identity network in $d$ dimensions. 
Even though this results in more neurons, it reduces the parameter count.
The following estimate is a consequence of \cref{lemma_parameters_of_id_conc}.

\begin{corollary}
\label{cor_parameters_of_id_conc}
Assume $a \in C(\R,\R)$ satisfies the $\idconstANN$-identity requirement 
for a number $\idconstANN \ge 2$ with $\IDANN{} \in \ANNs$ and denote $I_d = \paralANN{}{\IDANN{},\dots,\IDANN{}}$ for all $d \in \N$.
Let $\phi \in \ANNs$ and abbreviate $m = \max\{\inDimANN{\phi},\outDimANN{\phi}\}$.
Then
\begin{equation*}
	\neuronsANN{\IDANN{\outDimANN{\phi}} \concANN \phi \concANN \IDANN{\inDimANN{\phi}}} \leq \tfrac{5}{6} \idconstANN m \neuronsANN{\phi} + \tfrac{29}{12} \idconstANN^2 m^2.
\end{equation*}
\end{corollary}

This will be used in our proofs to estimate the number of parameters of 
\[
I_{l^{\phi_2}_{{\cal D}(\phi_2)}} \circ \phi_2 \circ I_{l^{\phi_2}_0}
\circ I_{l^{\phi_1}_{{\cal D}(\phi_1)}} \circ \phi_1 \circ I_{l^{\phi_1}_0}.
\]

\section{Catalog Networks}
\label{section_CNs}

\begin{figure}
	\centering
	\includegraphics{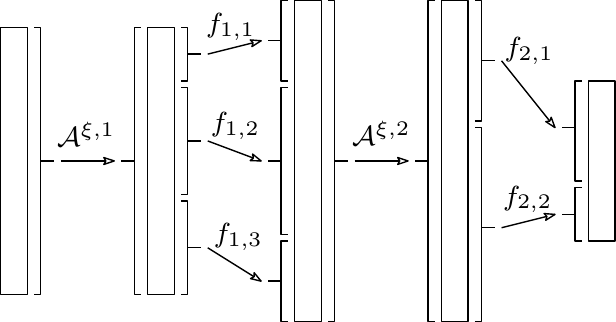}
	\caption{Realization of an example catalog network.}
	\label{fig_catalog_network}
\end{figure}

In this section, we generalize the concept of a neural network by allowing the activation functions 
to change from one layer to the next as long as they belong to a predefined catalog 
$\Catalog \subseteq \bigcup_{m,n \in \N} C(\R^m,\R^n)$. 
We denote the dimension of the domain of a function $f \in \bigcup_{m,n \in \N} C(\R^m,\R^n)$ by 
$\CatDimD{f}$ and the dimension of its target space by $\CatDimT{f}$, so that 
$f \in C(\R^{\CatDimD{f}},\R^{\CatDimT{f}})$. For a catalog $\Catalog$ and numbers $D \in \N$, $l_0,\dots,l_{2D} \in \N$, 
we define $\CNfixed$ as
\begin{equation*}
\begin{split}
	\prod_{k=1}^{D} \bigg[ \R^{l_{2k-1} \times l_{2k-2}} \times \R^{l_{2k-1}}
	\times \bigcup_{n \in \N} \Big\{ (f_1,\dots,f_n) \in \Catalog^n \colon \\
	\sum_{j=1}^{n} \CatDimD{f_j} = l_{2k-1} \text{ and } \sum_{j=1}^{n} \CatDimT{f_j} = l_{2k} \Big\} \bigg].
\end{split}
\end{equation*}
The set of all catalog networks corresponding to $\Catalog$ is given by 
\begin{equation*}
	\CNarb = \bigcup_{D \in \N} \bigcup_{l_0,\dots,l_{2D} \in \N} \CNfixed.
\end{equation*}
An element $\CNelt \in \CNfixed$ is of the form $\CNelt = \CNelementfixed$. For each $k \in \{1, \dots, D\}$,
we let $\CNaff{k}{\CNelt} \in C(\R^{l_{2k-2}},\R^{l_{2k-1}})$ be the affine function $x \mapsto V_k x+b_k$.
By $\CNnonlin{k}{\CNelt} \in C(\R^{l_{2k-1}},\R^{l_{2k}})$, we denote the function mapping 
$x \in \R^{l_{2k-1}}$ to
\begin{equation*}
\begin{split}
	\CNnonlin{k}{\CNelt}(x) = \big[ &f_{k,1} \big( x_{1}, \ldots, x_{\CatDimD{f_{k,1}}} \big), \\
	&f_{k,2} \big( x_{\CatDimD{f_{k,1}} + 1}, \ldots, x_{\CatDimD{f_{k,1}} + \CatDimD{f_{k,2}}} \big), \dots, \\
	&f_{k,n_k} \big( x_{\CatDimD{f_{k,1}} + \cdots + \CatDimD{f_{k,n_k-1}} + 1}, \dots, \\
	&\qquad \qquad \qquad \qquad x_{\CatDimD{f_{k,1}} + \cdots + \CatDimD{f_{k,n_k}}} \big) \big],
\end{split}
\end{equation*}
that is, we apply $f_{k,1}$ to the first $\CatDimD{f_{k,1}}$ entries of $x$, 
$f_{k,2}$ to the next $\CatDimD{f_{k,2}}$ entries and so on; see Figure \ref{fig_catalog_network}.
This is well-defined due to the sum conditions in the definition of $\CNfixed$.
The overall realization function $\CNfct{\CNelt} \in C(\R^{l_0},\R^{l_{2D}})$ of the catalog network 
$\CNelt$ is 
\begin{equation*}
	\CNfct{\CNelt} = \CNnonlin{D}{\CNelt} \circ \CNaff{D}{\CNelt} \circ \cdots \circ \CNnonlin{1}{\CNelt} \circ \CNaff{1}{\CNelt}.
\end{equation*}
We define the depth of $\CNelt$ as $\CNdepth{\CNelt} = D$.
Its input dimension is $\CNinput{\CNelt} = l_0$, its output dimension $\CNoutput{\CNelt} = l_{2D}$,
and its maximal width $\CNwidth{\CNelt} = \max\{l_0,\dots,l_{2D}\}$.

Our goal is to show that catalog networks can efficiently be approximated with neural networks
with respect to some weight function, by which we mean any function $\weight \colon [0,\infty) \rightarrow (0,\infty)$.

\begin{definition}
\label{def_weights_controlled_decay}
	We say the decay of a weight function $\weight$ is controlled by $(\OoGone,\OoGtwo) \in [1,\infty) \times [0,\infty)$ if
\begin{equation*}
\OoGone r^{\OoGtwo} \weight(r\max\{x,1\}) \ge \weight(x) 
\end{equation*}
for all $x \in [0,\infty)$ and $r \in [1,\infty)$.
\end{definition}

Controlled decay is a general concept applicable to different types of weight functions.
The inequality in Definition \ref{def_weights_controlled_decay} is exactly what is needed in the proofs of our
results. Useful weight functions are constants and functions of the form $(1+x^q)^{-1}$ or 
$(\max\{1,x^q\})^{-1}$ for some $q \in (0,\infty)$. Constant weight functions have decay controlled by $(1,0)$.
The functions $(1+x^q)^{-1}$ and $(\max\{1,x^q\})^{-1}$ are covered by the following result.

\begin{lemma}
\label{ex_weights_OoG}
Let $\delta \in (0,\infty)$ and consider a non-decreasing function $f \colon [0,\infty) \rightarrow (0,\infty)$.
Moreover, let $g \colon [0,\infty) \rightarrow [0,\infty)$ be of the form 
$x \mapsto \ssum{j=0}{q} a_j x^{b_j}$ for $q \in \N_0 = \N \cup \{0\}$ and $a_0,b_0,\dots,a_q,b_q \in [0,\infty)$.
Then the decay of the weight function $\weight(x) = f(x) (\max\{g(x),\delta\})^{-1}$ is controlled by 
$(\max\{g(1)/\delta,1\},\max\{b_0,\dots,b_q\})$.
\end{lemma}

Our main interest is in catalogs of functions that are well approximable with neural networks.
For the proofs of our main results to work we need the 
approximations to be Lipschitz continuous with a Lipschitz constant 
independent of the accuracy. To make this precise, we denote
the Euclidean norm by $\mednorm{\cdot}$.

\begin{definition}
\label{def_Cost}
Consider an activation function $a \in C(\R,\R)$ and a weight function $\weight$.
Fix constants $\Lip \in [0,\infty)$ and $\epsilon \in (0,1]$. Given a function $f \in \bigcup_{m,n \in \N} C(\R^m,\R^n)$ and 
a set $\approxSet{} \subseteq \R^{\CatDimD{f}}$, we define the approximation cost $\CostApprox{a}{\weight}{f}{\approxSet{}}{\Lip}{\epsilon}$ as the infimum of the set
\begin{equation*}
	\left\{ \neuronsANN{\phi} \in \N \colon
	\begin{array}{cc}
         \phi \in \ANNs \text{ with } \fctANN{\phi} \in C(\R^{\CatDimD{f}},\R^{\CatDimT{f}}) \\
         \text{s.t. } \fctANN{\phi} \text{ is } L\text{-Lipschitz on } \R^{\CatDimD{f}} \text{ and}\\
         \sup_{x \in \approxSet{}} \weight(\mednorm{x}) \mednorm{f(x) -  \fctANN{\phi}(x)} \leq \varepsilon
    \end{array}
\right\},
\end{equation*}
where, as usual, $\inf(\emptyset)$ is understood as $\infty$.
\end{definition}

The next definition specifies the class of catalogs for which we will be able
to prove \cref{theorem_main} on the approximability of catalog networks.

\begin{definition}
\label{def_approx_catalog}
Let $a \in C(\R,\R)$, $\kappa = (\kappa_0,\kappa_1,\kappa_2,\kappa_3) \in [1,\infty)^2 \times [0,\infty)^2$, 
$\epsilon \in (0,1]$, and suppose $\weight$ is a weight function.
Consider a subset $\Catalog \subseteq \bigcup_{m,n \in \N} C(\R^m,\R^n)$ together with a family of sets 
$\approxSet{} = (\approxSet{f})_{f \in \Catalog}$ such that $\approxSet{f} \subseteq \R^{\CatDimD{f}}$ contains $0$
for all $f \in {\cal F}$ and a collection of Lipschitz constants $\Lip = (\Lip_f)_{f \in \Catalog} \subseteq [0,\infty)$.
Then we call $\Catalog$ an $\ApproxCat{\epsilon}{\kappa}$-approximable catalog
if $\sup_{f \in \Catalog} \mednorm{f(0)} \leq \kappa_0$ and
\begin{equation*}
\CostApprox{a}{\weight}{f}{\approxSet{f}}{\Lip_f}{\delta} \leq \kappa_1 \max\{\CatDimD{f}, \CatDimT{f}\}^{\kappa_2} \delta^{-\kappa_3}
\end{equation*}
for all $f \in \Catalog$ and $\delta \in (0,\epsilon]$.
\end{definition}

Note that if $\Catalog$ is $\ApproxCat{\epsilon}{\kappa}$-approximable, then every $f \in \Catalog$ 
must be $\Lip_f$-Lipschitz continuous on the set $\approxSet{f}$.
Indeed, the definition implies that for all $\delta \in (0,\epsilon]$ there exists $\phi_{\delta} \in \ANNs$ such that 
$\weight(\mednorm{x}) \mednorm{f(x) - \fctANN{\phi_{\delta}}(x)} \leq \delta$ and 
$\mednorm{\fctANN{\phi_{\delta}}(x) - \fctANN{\phi_{\delta}}(y)} \leq \Lip_f \mednorm{x-y}$ for all $x,y \in \approxSet{f}$. 
Hence, one obtains from the triangle inequality that
\begin{equation*}
	\mednorm{f(x) - f(y)} \leq \frac{\delta}{\weight(\mednorm{x})} + \Lip_f \mednorm{x-y} + \frac{\delta}{\weight(\mednorm{y})}
\end{equation*}
for all $x,y \in \approxSet{f}$ and $\delta > 0$, which shows that $f$ is $L_f$-Lipschitz on $\approxSet{f}$.

If ${\cal F}$ is a catalog approximable on sets $\approxSet{} = (\approxSet{f})_{f \in \Catalog}$ with
Lipschitz constants $\Lip = (\Lip_f)_{f \in \Catalog} \subseteq [0,\infty)$, we define for a catalog network
$\CNelt \in \CNfixed$ of the form $\CNelt = \CNelementfixed$,
\begin{equation*}
\CNrdlayer{k}{\CNelt} = \sprod{j=1}{n_k} \approxSet{f_{k,j}} \subseteq \sprod{j=1}{n_k} \R^{\CatDimD{f_{k,j}}} = \R^{l_{2k-1}}
\end{equation*}
and
\begin{equation*}
\CNlayerLip{k}{\CNelt} = \max_{j \in \{1,\dots,n_k\}} \Lip_{f_{k,j}}
\end{equation*}
for all $k \in \{1,\dots,\CNdepth{\CNelt}\}$. Then the following holds.

\begin{lemma}
\label{lemma_main_layer_Lipschitz}
Let $\CNelt \in \CNarb$ be a catalog network based on an
$\ApproxCat{\epsilon}{\kappa}$-approximable catalog $\Catalog$.
Then 
\begin{equation*}
	\mednorm{\CNnonlin{k}{\CNelt}(x) - \CNnonlin{k}{\CNelt}(y)} \leq \CNlayerLip{k}{\CNelt}  \mednorm{x-y}
\end{equation*}
for all $k \in \{1,\dots,\CNdepth{\CNelt}\}$ and $x,y \in \CNrdlayer{k}{\CNelt}$.
\end{lemma}


\section{Examples of Approximable Catalogs}
\label{section_ex_approx_cat}

In this section, we provide different examples of approximable catalogs that will be used in \cref{section_applications} to 
show that various high-dimensional functions admit neural network approximations without the 
curse of dimensionality. Our catalogs are based on one-dimensional Lipschitz functions, the maximum function, 
the square, the product and the decreasing exponential function. They will be collected in \cref{ex_Lip_catalog,ex_Lip_max_catalog,ex_RBF_catalog,ex_Lip_prod_catalog}.

First, consider a $\Lipconst$-Lipschitz function 
$f \colon \R \rightarrow \R$ for a constant $\Lipconst \in [0,\infty)$.
For any given $r \in (0,\infty)$, $f$ can be approximated on $[-r,r]$ with a piece-wise linear function 
supported on $N+1$ equidistributed points with accuracy $\Lipconst r/N$.
Such a piece-wise linear function can be realized with a ReLU network $\phi_N$ with one hidden layer
and $N$ hidden neurons. This results in $\neuronsANN{\phi_N} = 3N+1$, from which it follows that
\begin{equation*}
	\CostApprox{\ReLU}{1}{f}{[-r,r]}{\Lipconst}{\epsilon} \leq \neuronsANN{\phi_{\lceil \Lipconst r\epsilon^{-1} \rceil}} \leq 3\Lipconst r\epsilon^{-1} + 4.
\end{equation*}
Alternatively, one can approximate $f$ on the entire real line with respect to a 
weight function of the form $\weight_q(x) = (1+x^q)^{-1}$ for some $q \in (1,\infty)$.
Then 
\begin{equation*}
	\CostApprox{\ReLU}{\weight_q}{f}{\R}{\Lipconst}{\epsilon} \leq 2^{\nicefrac{1}{q-1}} 3 (\Lipconst \epsilon^{-1})^{\nicefrac{q}{q-1}} + 4,
\end{equation*}
the proof of which is a variant of \cite[Corollary\,3.13]{HutzJenKruseNguyen2020}.
Indeed, set $r = (2\Lipconst \epsilon^{-1})^{\nicefrac{1}{(q-1)}}$ and $N = \lceil \Lipconst r\epsilon^{-1}\rceil$.
Using $\phi_N$ as above, we have $|f(x)-\phi_N(x)| \leq \epsilon$ for all $x \in [-r,r]$ and $|f(x)-\phi_N(x)| \leq 2\Lipconst |x|$ for all $x \in \R \backslash [-r,r]$. The choice of $r$ then ensures that $\weight_q(|x|)|f(x)-\phi_N(x)| \leq \epsilon$ for all $x \in \R$.
In the notation of approximable catalogs, we can summarize as follows.

\begin{example}
\label{ex_Lip_catalog}
Let $r \in (0,\infty)$ and consider the weight function $\weight_q(x) = (1+x^q)^{-1}$
for a $q \in (1,\infty)$. For $\Lipconst \in [1,\infty)$, introduce the $\Lipconst$-Lipschitz catalog
\begin{equation*}
\begin{split}
	\Lipcatalog = \big\{ f \in C(\R,\R) \colon &f \text{ is } \Lipconst \text{-Lipschitz} 
	\text{ and } |f(0)| \leq \Lipconst \big\}.
\end{split}
\end{equation*}
Set $\Lip_{\ID{}} = 1$ and $\Lip_f = \Lipconst$ for $f \in \Lipcatalog \backslash \{\ID{}\}$.
If we define approximation sets by $\approxSet{\ID{}} = \R$ and
\begin{enumerate}[\rm 1)]
\item \label{ExLip(i)}
$\approxSet{f} = [-r,r]$ for all $f \in \Lipcatalog \backslash \{\ID{}\}$, then $\Lipcatalog$ is a $\ApproxExCat{1}{\kappa}{1}$-approximable catalog for $\kappa = (\Lipconst,3\Lipconst r + 4,0,1)$

\item \label{ExLip(ii)}
$\approxSet{f} = \R$ for all $f \in \Lipcatalog \backslash \{\ID{}\}$, then $\Lipcatalog$ is a $\ApproxExCat{1}{\kappa}{\weight_q}$-approximable catalog for\footnote{Here we use that $4 \leq 2 (2\Lipconst \epsilon^{-1})^{\nicefrac{q}{(q-1)}}$ since $\Lipconst \geq 1$.} $\kappa = (\Lipconst,5 (2\Lipconst)^{\nicefrac{q}{(q-1)}},0, \nicefrac{q}{q-1})$.
\end{enumerate}
\end{example}

Let us now turn to the maximum functions $\maxfct{d} \colon \R^d \rightarrow \R, ~x \mapsto \max\{x_1,\dots,x_d\}$, $d \in \N$.
They admit an exact representation with ReLU networks.
Indeed, $\maxfct{1}$ is simply the identity and 
$\maxfct{2}$ is the ReLU-realization of
\begin{equation*}
	\phi_2 = \left[ \left(
	\begin{bmatrix}
		1 & -1 \\
		0 & 1 \\
		0 & -1
	\end{bmatrix},
	\begin{bmatrix}
		0 \\
		0 \\
		0
	\end{bmatrix} \right), \big(
	\begin{bmatrix}
		1 & 1 & -1
	\end{bmatrix},
	0 \big) \right]\!.
\end{equation*}
If $\IDANN{} \in \ANNs$ is a skeleton for which ReLU satisfies the 
2-identity requirement and we define $\IDANN{d} = \paralANN{}{\IDANN{},\dots,\IDANN{}}$, $d \in \N$, 
then it easily follows by induction that $\maxfct{d}$, $d \in \N_{\geq 3}$, is the ReLU-realization of 
$\phi_d = \phi_{d-1} \concANN \paralANN{}{\phi_2,\IDANN{d-2}}$, whose architecture is
$(d,2d-1,2d-3,\dots,3,1)$. From this, we obtain $\neuronsANN{\phi_d} = \frac{1}{3}(4d^3 + 3d^2 - 4d + 3) \leq 2d^3$.
In other words, for all $d \in \N$ and any weight function $\weight$,
\begin{equation*}
\CostApprox{\ReLU}{\weight}{\maxfct{d}}{\R^d}{1}{\epsilon} \leq 2d^3.
\end{equation*}
Adding the maximum functions to the Lipschitz catalog, we obtain the following.

\begin{example}
\label{ex_Lip_max_catalog}
Adopt the setting of \cref{ex_Lip_catalog} and define the $\Lipconst$-Lipschitz-maximum catalog 
$\Lipmaxcatalog = \Lipcatalog \cup \{\maxfct{d} \colon d \in \N\}$. Add the approximation set $\approxSet{\maxfct{d}} = \R^d$
and the Lipschitz constant $\Lip_{\maxfct{d}} = 1$ for all $d \in \N$.
Then $\Lipmaxcatalog$ is
\begin{enumerate}[\rm 1)] 
\item \label{max(i)} a $\ApproxExCat{1}{\kappa}{1}$-approximable catalog for $\kappa = (\Lipconst,3\Lipconst r + 4,3,1)$
and $\approxSet{}$ as in \cref{ex_Lip_catalog}.\eqref{ExLip(i)}.

\item \label{max(ii)} a $\ApproxExCat{1}{\kappa}{\weight_q}$-approximable catalog for $\kappa =
(\Lipconst,5 (2\Lipconst)^{\nicefrac{q}{(q-1)}},3,\nicefrac{q}{q-1})$ and $\approxSet{}$ as in 
\cref{ex_Lip_catalog}.\eqref{ExLip(ii)}.
\end{enumerate}
\end{example}

Next, we study the approximability of the square function $\sqfct \colon \R \rightarrow \R, ~x \mapsto x^2$ .
It has been shown by different authors that it can be approximated with accuracy $\epsilon > 0$ on the unit interval 
by the ReLU-realization of a skeleton $\phi_{\epsilon} \in {\cal N}$ 
satisfying $\neuronsANN{\phi_{\epsilon}} = \mathcal{O}(\log_2(\epsilon^{-1}))$; see \cite{Yarotsky2017,SchwabZech2019,GrohsHorJenZimm2019,GrohsPerekElbBolc2019}.
A precise estimate of the required number of parameters is given in Proposition 3.3 of
\cite{GrohsHorJenZimm2019}. In our language it can be stated as
\begin{equation*}
\CostApprox{\ReLU}{1}{\sqfct}{[0,1]}{2}{\epsilon} \leq \max\{13,10\log_2(\epsilon^{-1})-7\}.
\end{equation*}
Moreover, the neural network $\fctReLUANN{\phi_{\epsilon}}$ achieving this cost is 
2-Lipschitz and satisfies $\fctReLUANN{\phi_{\epsilon}} = \ReLU$ on $\R \backslash [0,1]$.
Using a mirroring and scaling argument, we can deduce the following estimate for 
approximating the square function on the interval $[-r,r]$ for any $r \in (0,\infty)$.

\begin{lemma}
\label{lemma_sq_mirror}
For all $r \in (0,\infty)$ and $\epsilon \in (0,1]$, there exists a skeleton
$\psi_{r,\epsilon} \in \ANNs$ such that $\fctReLUANN{\psi_{r,\epsilon}} \in C(\R, \R)$ is $2r$-Lipschitz,
$\sup_{x \in [-r,r]} |\fctReLUANN{\psi_{r,\epsilon}}(x) - x^2| \leq \epsilon$, 
$\fctReLUANN{\psi_{r,\epsilon}}(x) = r|x|$ for all $x \in \R \backslash [-r,r]$ and
\begin{equation*}
	\neuronsANN{\psi_{r,\epsilon}} \leq \max\{52,80\log_2(r) + 40\log_2(\epsilon^{-1}) - 28\}.
\end{equation*}
\end{lemma}

More concisely, for all $r \in [2,\infty)$ and $\epsilon \in (0,1]$, the statement of Lemma \ref{lemma_sq_mirror} can be written as
\begin{equation*}
\begin{split}
	\CostApprox{\ReLU}{1}{\sqfct}{[-r,r]}{2r}{\epsilon} &\leq \\
	80\log_2(r) &+ 40\log_2(\epsilon^{-1}) - 28.
\end{split}
\end{equation*}

Now, let us take a closer look at the decreasing exponential function 
$\expfct \colon \R \rightarrow \R, ~x \mapsto e^{-x}$. Its restriction to $[0,\infty)$ is
covered by the general approximation result for Lipschitz functions. But exploiting its 
exponential decrease, we can obtain better estimates. More precisely, $\expfct$
can be approximated to a given accuracy $\epsilon \in (0,1]$ uniformly on $[0,\infty)$ with a piece-wise linear 
interpolation supported on the $\lfloor \epsilon^{-1} \rfloor$ points $-\log(n\epsilon)$, $n \in \{1,\dots,\lfloor\epsilon^{-1}\rfloor \}$
which is constant on $\R \backslash [-\log(\lfloor\epsilon^{-1}\rfloor \epsilon),-\log(\epsilon)]$.
Realizing this piece-wise linear function with a ReLU network with one hidden layer yields
\begin{equation*}
\CostApprox{\ReLU}{1}{\expfct}{[0,\infty)}{1}{\epsilon} \leq 3 \lfloor \epsilon^{-1} \rfloor + 1 \leq 4 \epsilon^{-1}.
\end{equation*}
Together with $\ID{}$ and $\sqfct$, $\expfct$ gives rise to the following catalog, which we will use to approximate
generalized Gaussian radial basis function networks in \cref{section_applications}.

\begin{example}
\label{ex_RBF_catalog}
Let $r \in [5,\infty)$. Define the catalog $\GaussianRBFcat = \{\ID{},\expfct,\sqfct\}$ with approximation sets $\approxSet{\ID{}} = \R$, 
$\approxSet{\expfct} = [0,\infty)$, $\approxSet{\sqfct} = [-r,r]$ and Lipschitz constants $\Lip_{\ID{}} = \Lip_{\expfct} = 1$, 
$\Lip_{\sqfct} = 2r$. Then $\GaussianRBFcat$ is\footnote{That $r \geq 5$ and $\epsilon \leq r^{-3}$ shows $80\log_2(r) + 40\log_2(\epsilon^{-1})-28 \leq 4\epsilon^{-1}$.} a $\ApproxExCat{r^{-3}}{(1,4,0,1)}{1}$-approximable catalog.
\end{example}

Using the identity $xy = \frac{1}{4}((x+y)^2 - (x-y)^2)$, we can also 
estimate the approximation rate of the product function $\prodfct \colon \R^2 \rightarrow \R, ~(x,y) \mapsto xy$.
This trick has already been used before by, e.g., \cite{Yarotsky2017,LinTegRol2017}. 
We still provide a proof of the following proposition since the results in the existing literature do not specify the Lipschitz constant.
Our proofs of both, \cref{lemma_sq_mirror} and \cref{prop_product_approx}, follow the reasoning of Section 3 in \cite{GrohsHorJenZimm2019}.

\begin{proposition}
\label{prop_product_approx}
For all $r \in (0,\infty)$ and $\epsilon \in (0,\frac{1}{2}]$, one has
\begin{equation*}
\begin{split}
	\CostApprox{\ReLU}{1}{\prodfct}{[-r,r]^2}{\sqrt{8}r}{\epsilon} &\leq \\
	\max\{208,320\log_2(r) &+ 160\log_2(\epsilon^{-1}) + 48\}.
\end{split}
\end{equation*}
\end{proposition}

The following is our last example of an approximable catalog. 

\begin{example}
\label{ex_Lip_prod_catalog}
Take $\Lipcatalog$ from \cref{ex_Lip_catalog}, let $R \in (0,\infty)$, and define the 
$\Lipconst$-Lipschitz-product catalog $\Lipprodcatalog = \Lipcatalog \cup \{\prodfct\}$.
The approximation sets and Lipschitz constants are defined as in
\eqref{ExLip(i)} of \cref{ex_Lip_catalog} for $\Lipcatalog$ and 
$\approxSet{\prodfct} = [-R,R]^2$, $\Lip_{\prodfct} = \sqrt{8}R$.
Then $\Lipprodcatalog$ is a $\ApproxExCat{\delta}{(\Lipconst,M,0,1)}{1}$-approximable catalog
for\footnote{This specific choice of $\delta$ and $M$ ensures that $\max\{208,320\log_2(R) + 160\log_2(\epsilon^{-1}) + 48\} 
\leq M \epsilon^{-1}$ for all $\epsilon \in (0,\delta]$.} $\delta = \min\{\nicefrac{1}{2},\nicefrac{R^2}{2}\}$ and 
$M = \max\{3\Lipconst r + 4,105R^2\}$.
\end{example}


\section{Approximation Results}
\label{section_approx_results}

In this section, we state the first of our main results, \cref{theorem_main},
on the approximability of catalog networks with neural networks and explore 
the special case of ReLU activation in \cref{cor_ReLU,cor_ReLU_weight_one}. 
The next lemma is crucial for the proof of \cref{theorem_main}. It establishes the approximability 
of the functions $\CNnonlin{k}{\CNelt}$, $k \in \{1,\dots,\CNdepth{\CNelt}\}$, 
in a catalog network $\CNelt \in \CNarb$. Since $\CNnonlin{k}{\CNelt}$ is composed of
functions $f_{k,1}, \dots, f_{k,n_k}$ from the catalog $\Catalog$, it can be approximated by approximating 
$f_{k,1}, \dots, f_{k,n_k}$ with neural networks and then parallelizing them as in Figure \ref{fig_parallelization_novel}.
\cref{prop_ANN_paral_NEW} allows us to keep track of the resulting number of parameters.

\begin{lemma}
\label{lemma_main_layer_approx}
Assume $a \in C(\R,\R)$ satisfies the $\idconstANN$-identity requirement for some $\idconstANN \geq 2$.
Let $\Catalog$ be an $\ApproxCat{\epsilon}{\kappa}$-approximable catalog for a non-increasing weight function $\weight$,
and consider a catalog network $\CNelt \in \CNfixed$ for some $D \in \N$ and $l_0,\dots,l_{2D} \in \N$.
Then for all $k \in \{1,\dots,D\}$ and $\delta \in (0,\epsilon]$, there exists a skeleton
$\phi \in \ANNs$ with $a$-realization $\fctANN{\phi} \in C(\R^{l_{2k-1}},\R^{l_{2k}})$ such that 
\begin{enumerate}[\rm 1)]
\item $\sup_{x \in \CNrdlayer{k}{\CNelt}} \weight(\mednorm{x}) \mednorm{\CNnonlin{k}{\CNelt}(x) - \fctANN{\phi}(x)} \leq \delta$,

\item $\fctANN{\phi}$ is $\CNlayerLip{k}{\CNelt}$-Lipschitz continuous on $\R^{l_{2k-1}}$ and

\item\label{lemma_main_layer_approx_item} $\neuronsANN{\phi} \leq \frac{11}{16} \kappa_1 
\idconstANN^2 \max\{l_{2k-1},l_{2k}\}^{\kappa_2 + \kappa_3/2 + 5} \delta^{-\kappa_3}$.
\end{enumerate}
If, in addition, $\CatDimD{f} \leq d$ and $\CatDimT{f} \leq d$ for some $d \in \N$ and all $f \in \Catalog$, then one also has
\begin{enumerate}[\rm 1)]
\setcounter{enumi}{3}
\item\label{lemma_main_layer_approx_item_prime} $\neuronsANN{\phi} \leq \frac{11}{16} 
\kappa_1 \idconstANN^2 d^2 \max\{l_{2k-1},l_{2k}\}^{\kappa_2 + \kappa_3/2 + 3} \delta^{-\kappa_3}$.
\end{enumerate}
\end{lemma}

Before we can formulate \cref{theorem_main}, we have to introduce a few more concepts. Let $\Catalog$ be a catalog that is 
approximable on sets $\approxSet{} = (\approxSet{f})_{f \in \Catalog}$ 
with Lipschitz constants $\Lip = (\Lip_f)_{f \in \Catalog} \subseteq [0,\infty)$.
Then, for any catalog network $\CNelt = \CNelementfixed \in \CNarb$, we define
\begin{equation*}
\begin{split}
	&\CNredom{\CNelt} = \bigg\{ x \in \R^{\CNinput{\CNelt}} \colon \text{for all } k \in \{1,\dots,\CNdepth{\CNelt}\} \colon \\
	&\big( \CNaff{k}{\CNelt} \circ \CNnonlin{k-1}{\CNelt} \circ \cdots \circ\CNaff{2}{\CNelt} \circ \CNnonlin{1}{\CNelt} \circ \CNaff{1}{\CNelt} \big) (x) \in \CNrdlayer{k}{\CNelt} \bigg\}
\end{split}
\end{equation*}
and 
\begin{equation*}
	\CNfullLip{\CNelt} = \prod_{k=1}^{\CNdepth{\CNelt}} \CNlayerLip{k}{\CNelt} \opnorm{V_k},
\end{equation*}
where $\opnorm{\cdot}$ denotes the operator norm when applied to matrices.
The set $\CNredom{\CNelt}$ describes where we will be able to approximate the catalog network $\xi$.
It takes into account that each each layer function $\CNnonlin{k}{\CNelt}$ can only be approximated on the set $\CNrdlayer{k}{\CNelt}$.
The number $\CNfullLip{\CNelt}$ represents the worst-case Lipschitz constant of the catalog network.

To estimate the approximation error, we need two more quantities.
The first one is
\begin{equation*}
	\CNtransl{\CNelt} = \max\big\{1, \mednorm{\CNaff{1}{\CNelt}(0)}, \dots, \mednorm{\CNaff{\CNdepth{\CNelt}}{\CNelt}(0)} \big\},
\end{equation*}
which simply measures the maximal norm of the inhomogeneous parts of the affine transformations 
(capped from below by 1). When using weight functions of the type $w_q(x) = (1+x^q)^{-1}$, 
functions in the catalog are approximated better close to the origin. The quantity $\CNtransl{\CNelt}$  
together with the $\kappa_0$-boundedness of the catalog in the origin
will be used to control how far away one is from the region where one has the best approximation.
However, this becomes irrelevant for constant weight functions, as can be seen in 
\cref{cor_ReLU_weight_one} below.

The last quantity we need is $\CNtechLip{\CNelt}$, defined as the maximum of 1 and
\begin{equation*}
	\max_{k \in \{0,\dots,D-1\}} \max\{1, \CNlayerLip{k}{\CNelt}, \CNlayerLip{D}{\CNelt}\} \opnorm{V_D} \prod_{j=k+1}^{D-1} \CNlayerLip{j}{\CNelt} \opnorm{V_j}
\end{equation*}
where we abbreviate $D = \CNdepth{\CNelt}$ and use the convention $\CNlayerLip{0}{\CNelt} = 0$.
This combines the Lipschitz constants of the affine and nonlinear functions appearing in the different layers of the 
catalog network $\xi$.

\begin{theorem}
\label{theorem_main}
Suppose $a \in C(\R,\R)$ fulfills the $\idconstANN$-identity requirement for some number $\idconstANN \ge 2$
and let $\weight$ be a non-increasing weight function whose decay is controlled by $(\OoGone,\OoGtwo)$ for some $\OoGone \in [1,\infty)$ and $\OoGtwo \in [0,\infty)$. Consider a catalog network $\CNelt \in \CNarb$ for
an $\ApproxCat{\epsilon}{\kappa}$-approximable catalog $\Catalog$.
Then there exists a skeleton $\phi \in \ANNs$ with $a$-realization
$\fctANN{\phi} \in C(\R^{\CNinput{\CNelt}},\R^{\CNoutput{\CNelt}})$ such that
\begin{enumerate}[\rm 1)]
\item $\sup_{x \in \CNredom{\CNelt}} \weight(\mednorm{x})\mednorm{\CNfct{\CNelt}(x) - \fctANN{\phi}(x)} \leq \epsilon$,

\item $\fctANN{\phi}$ is $\CNfullLip{\CNelt}$-Lipschitz continuous on $\R^{\CNinput{\CNelt}}$ and 

\item\label{theorem_main_item}
$\neuronsANN{\phi} \leq C \CNtransl{\CNelt}^{t-\kappa_3} \CNtechLip{\CNelt}^t \CNdepth{\CNelt}^{t+1} \CNwidth{\CNelt}^{\kappa_2 + t/2 + 6} \epsilon^{-\kappa_3}$ \newline for $t = \kappa_3(\OoGtwo + 1)$ and $C = \tfrac{81}{32} \idconstANN^3 (4 \kappa_0)^{t-\kappa_3} \kappa_1 \OoGone^{\kappa_3}$.
\end{enumerate}
\end{theorem}

The conclusion of \cref{theorem_main} could be written more concisely as
\begin{equation*}
\begin{split}
\CostApprox{a}{\weight}{\CNfct{\CNelt}}{&\,\CNredom{\CNelt}}{\CNfullLip{\CNelt}}{\epsilon} \\
&\leq C \CNtransl{\CNelt}^{t-\kappa_3} \CNtechLip{\CNelt}^t \CNdepth{\CNelt}^{t+1} 
\CNwidth{\CNelt}^{\kappa_2 + t/2 + 6} \epsilon^{-\kappa_3}.
\end{split}
\end{equation*}
We point out that the rate in the accuracy is $\mathcal{O}(\epsilon^{-\kappa_3})$, the same 
as for the underlying catalog $\Catalog$.

In the proof of \cref{theorem_main}, we combine the approximations of the functions $\CNnonlin{k}{\CNelt}$
obtained in \cref{lemma_main_layer_approx} with the affine maps $\CNaff{k}{\CNelt}$. When concatenating
different approximating networks, we interpose identity networks. This reduces the parameter count 
in worst-case scenarios but can lead to slightly looser estimates in certain other situations.

\begin{remark}
{\rm \label{rem_improved_exp_for_dim_one_thrm}
If $\CatDimD{f} \leq d$ and $\CatDimT{f} \leq d$ for some $d \in \N$ and all $f \in \Catalog$, we can use (\ref{lemma_main_layer_approx_item_prime}) instead of (\ref{lemma_main_layer_approx_item}) of
\cref{lemma_main_layer_approx} in the proof of \cref{theorem_main} to obtain the following modified version of the parameter bound in \cref{theorem_main}:
\begin{equation*}
\neuronsANN{\phi} \leq C d^2 \CNtransl{\CNelt}^{t-\kappa_3} \CNtechLip{\CNelt}^t \CNdepth{\CNelt}^{t+1} \CNwidth{\CNelt}^{\kappa_2 + t/2 + 4} \epsilon^{-\kappa_3}.
\end{equation*}
Since in many of the example catalogs of \cref{section_ex_approx_cat}, the maximal input/output dimension is 1 or 2,
this will allow us to obtain better estimates in some of the applications in \cref{section_applications} below.
}
\end{remark}

\begin{remark}
\label{rem_norm}
{\rm
A careful inspection of the proof of Theorem \ref{theorem_main} shows that it does not only work for
the Euclidean norm but also, for instance, the sup-norm.}
\end{remark}

We know that the ReLU activation function satisfies the 2-identity requirement. \cref{theorem_main} 
recast for ReLU activation and the weight function $w_q(x) = (1+x^q)^{-1}$ reads as follows:

\begin{corollary}
\label{cor_ReLU}
Consider the weight function $\weight_q(x) = (1+x^q)^{-1}$ for some
$q \in (0,\infty)$. Let $\CNelt \in \CNarb$ be a catalog network for 
a $\ApproxExCat{\epsilon}{\kappa}{\weight_q}$-approximable catalog $\Catalog$.
Then there exists a skeleton $\phi \in \ANNs$ with ReLU-realization 
$\fctReLUANN{\phi} \in C(\R^{\CNinput{\CNelt}},\R^{\CNoutput{\CNelt}})$ such that
\begin{enumerate}[\rm 1)]
\item $\sup_{x \in \CNredom{\CNelt}} (1+\mednorm{x}^q)^{-1}\mednorm{\CNfct{\CNelt}(x) - \fctReLUANN{\phi}(x)} \leq \epsilon$,

\item $\fctReLUANN{\phi}$ is $\CNfullLip{\CNelt}$-Lipschitz continuous on $\R^{\CNinput{\CNelt}}$ and 

\item 
$\neuronsANN{\phi} \leq C \CNtransl{\CNelt}^{t-\kappa_3} \CNtechLip{\CNelt}^t \CNdepth{\CNelt}^{t+1} \CNwidth{\CNelt}^{\kappa_2 + t/2 + 6} \epsilon^{-\kappa_3}$ \newline
for $t = \kappa_3(q + 1)$ and $C = \tfrac{81}{4} 2^{2t-\kappa_3} \kappa_0^{t-\kappa_3} \kappa_1$.
\end{enumerate}
\end{corollary}

For the weight function $w \equiv 1$, the parameter estimate in Theorem  \ref{theorem_main} simplifies considerably.
This is because the decay of $w \equiv 1$ is controlled by $(1,0)$, which makes the translation size and the bound of the 
catalog in the origin irrelevant.

\begin{corollary}
\label{cor_ReLU_weight_one}
Let $\CNelt \in \CNarb$ be a catalog network for a $\ApproxExCat{\epsilon}{\kappa}{1}$-approximable catalog $\Catalog$.
Then there exists a skeleton $\phi \in \ANNs$ with ReLU-realization 
$\fctReLUANN{\phi} \in C(\R^{\CNinput{\CNelt}},\R^{\CNoutput{\CNelt}})$ such that 
\begin{enumerate}[\rm 1)]
\item $\sup_{x \in \CNredom{\CNelt}} \mednorm{\CNfct{\CNelt}(x) - \fctReLUANN{\phi}(x)} \leq \epsilon$,

\item $\fctReLUANN{\phi}$ is $\CNfullLip{\CNelt}$-Lipschitz continuous on $\R^{\CNinput{\CNelt}}$ and 

\item $\neuronsANN{\phi} \leq \tfrac{81}{4} \kappa_1 \CNtechLip{\CNelt}^{\kappa_3} \CNdepth{\CNelt}^{\kappa_3+1} \CNwidth{\CNelt}^{\kappa_2 + \kappa_3/2 + 6} \epsilon^{-\kappa_3}$.
\end{enumerate}
\end{corollary}


\section{Log-approximable Catalogs}
\label{section_log_modification}

In this section, we modify the way we measure the approximation cost and derive corresponding 
approximation results.

\begin{definition}
\label{def_approx_catalog_log}
With the setup of \cref{def_approx_catalog} and $\epsilon \leq 1/2$,
a subset $\Catalog \subseteq \bigcup_{m,n \in \N} C(\R^m,\R^n)$ is called
$\ApproxCat{\epsilon}{\kappa}$-log-approximable
if $\sup_{f \in \Catalog} \mednorm{f(0)} \leq \kappa_0$ and
\begin{equation*}
	\CostApprox{a}{\weight}{f}{\approxSet{f}}{\Lip_f}{\delta} \leq \kappa_1\! \max\{\CatDimD{f}, \CatDimT{f}\}^{\!\kappa_2} [\log_2(\delta^{-1})]^{\kappa_3}
\end{equation*}
for all $f \in \Catalog$ and $\delta \in (0,\epsilon]$.
\end{definition}

This log-modification is designed for catalogs made of functions like the square or the product, which 
can be approximated with rate $\mathcal{O}(\log_2(\epsilon^{-1}))$, as we have seen in \cref{lemma_sq_mirror,prop_product_approx}.
Its usefulness will become apparent in \cref{ex_prod}, which is based on the following catalog.

\begin{example}
\label{ex_prod_catalog}
{\rm Let $\Lipprod = \{\ID{},\prodfct\}$ be the product catalog and fix $r \in [1,\infty)$ and $d \in \N$.
Consider the approximation sets $\approxSet{\ID{}} = \R$, $\approxSet{\prodfct} = [-r^d,r^d]^2$ and 
the Lipschitz constants $\Lip_{\ID{}} = 1$, $\Lip_{\prodfct} = \sqrt{8}r^d$.
Then $\Lipprod$ is a $\ApproxExCat{\min\{1/2,1/r\}}{(1,320d+208,0,1)}{1}$-log-approximable 
catalog, where $\kappa_1$ can also be chosen as $208$ instead of $320d+208$ if $r=1$.}
\end{example}

Using the notion of log-approximable catalogs, we can derive the following analogue of \cref{theorem_main}.

\begin{theorem}
\label{theorem_main_LOG}
Assume $a \in C(\R,\R)$ satisfies the $\idconstANN$-identity requirement for some number $\idconstANN \ge 2$,
and let $\weight$ be a non-increasing weight function whose decay is controlled by $(\OoGone,\OoGtwo)$ for some $\OoGone \in [1,\infty)$ and $\OoGtwo \in [0,\infty)$. Consider a catalog network $\CNelt \in \CNarb$ for an
$\ApproxCat{\epsilon}{\kappa}$-log-approximable catalog ${\cal F}$.
Then there exists a skeleton $\phi \in \ANNs$ with $a$-realization 
$\fctANN{\phi} \in C(\R^{\CNinput{\CNelt}},\R^{\CNoutput{\CNelt}})$ such that 
\begin{enumerate}[\rm 1)]
\item $\sup_{x \in \CNredom{\CNelt}} \weight(\mednorm{x})\mednorm{\CNfct{\CNelt}(x) - \fctANN{\phi}(x)} \leq \epsilon$,

\item $\fctANN{\phi}$ is $\CNfullLip{\CNelt}$-Lipschitz continuous on $\R^{\CNinput{\CNelt}}$ and 

\item 
$\neuronsANN{\phi} \leq \tfrac{81}{32} \idconstANN^3 \kappa_1 \CNdepth{\CNelt} \CNwidth{\CNelt}^{\kappa_2 + 6}$ \newline $\times \big[\!\log_2\!\big( \OoGone [4 \kappa_0 \CNtransl{\CNelt}]^{\OoGtwo} [\CNtechLip{\CNelt} \CNdepth{\CNelt} \sqrt{\CNwidth{\CNelt}}]^{\OoGtwo +1} \epsilon^{-1} \big)\big]^{\kappa_3}$.

\end{enumerate}
\end{theorem}

The following is the analogue of \cref{cor_ReLU_weight_one} for log-approximable catalogs.

\begin{corollary}
\label{cor_ReLU_weight_one_LOG}
Let $\CNelt \in \CNarb$ be a catalog network 
for a $\ApproxExCat{\epsilon}{\kappa}{1}$-log-approximable catalog $\Catalog$.
Then there exists a skeleton $\phi \in \ANNs$ with 
ReLU-realization $\fctReLUANN{\phi} \in C(\R^{\CNinput{\CNelt}},\R^{\CNoutput{\CNelt}})$ such that 
\begin{enumerate}[\rm 1)]
\item $\sup_{x \in \CNredom{\CNelt}} \mednorm{\CNfct{\CNelt}(x) - \fctReLUANN{\phi}(x)} \leq \epsilon$,

\item $\fctReLUANN{\phi}$ is $\CNfullLip{\CNelt}$-Lipschitz continuous on $\R^{\CNinput{\CNelt}}$ and 

\item $\neuronsANN{\phi} \leq \tfrac{81}{4} \kappa_1 \CNdepth{\CNelt} \CNwidth{\CNelt}^{\kappa_2 + 6} \big[\!\log_2\!\big( \CNtechLip{\CNelt} \CNdepth{\CNelt} \sqrt{\CNwidth{\CNelt}} \epsilon^{-1} \big)\big]^{\kappa_3}$.
\end{enumerate}
\end{corollary}


\section{Overcoming the Curse of Dimensionality}
\label{section_applications}

In this section, we apply the theory of catalog networks to show that different high-dimensional functions 
admit a ReLU neural network approximation without the curse of dimensionality.
We use the catalogs introduced in \cref{section_ex_approx_cat,section_log_modification} to construct
families of functions indexed by the dimension of their domain that are of the same form for each dimension.
The results in this section are proved by finding catalog network representations of the high-dimensional target functions,
so that one of the general approximation results of Sections \ref{section_approx_results} or \ref{section_log_modification} 
can be applied. The mere approximability of these functions with ReLU networks follows from classical 
universal approximation results such as 
\cite{LeshnoLinPinkusScho1993}. But the quantitative estimates on the number of parameters in terms of
the dimension $d$ and the accuracy $\epsilon$ are new.

Our general results cover a wide range of interesting examples, but it is possible that, for some of them, 
the estimates could be improved by using their special structure.

\begin{proposition}
\label{ex_Lipschitz_sum_bounded}
Fix $\Lipconst,r \in [1,\infty)$, and let $f_i \colon \R \rightarrow \R$, $i \in \N$, be 
$\Lipconst$-Lipschitz continuous with $|f_i(0)| \leq \Lipconst$.
Define $g_d \colon \R^d \rightarrow \R$ by $g_d(x) = \ssum{i=1}{d} f_i(x_i)$.
Then for all $d \in \N$ and $\epsilon \in (0,1]$, there exists $\phi \in \ANNs$ with $\fctReLUANN{\phi} \in C(\R^d,\R)$ such that
\begin{enumerate}[\rm 1)]
\item $\sup_{x \in [-r,r]^d} |g_d(x) - \fctReLUANN{\phi}(x)| \leq \epsilon$,

\item $\fctReLUANN{\phi}$ is $\sqrt{d}\Lipconst$-Lipschitz continuous on $\R^d$ and

\item $\neuronsANN{\phi} \leq \tfrac{4}{7} 10^3 \Lipconst^2 r d^5 \epsilon^{-1}$.
\end{enumerate}
\end{proposition}

The following is a generalization, whose proof only requires a slight adjustment.

\begin{proposition}
\label{ex_Lipschitz_sum_bounded_v2}
Fix $\Lipconst,r \in [1,\infty)$, and let $f_i \colon \R \rightarrow \R$, $i \in \N_0$, be 
$\Lipconst$-Lipschitz continuous with $|f_i(0)| \leq \Lipconst$.
Define $g_d \colon \R^d \rightarrow \R$ by $g_d(x) = f_0(\ssum{i=1}{d} f_i(x_i))$.
Then for all $d \in \N$ and $\epsilon \in (0,1]$, there exists $\phi \in \ANNs$ with $\fctReLUANN{\phi} \in C(\R^d,\R)$ such that
\begin{enumerate}[\rm 1)]
\item $\sup_{x \in [-r,r]^d} |g_d(x) - \fctReLUANN{\phi}(x)| \leq \epsilon$,

\item $\fctReLUANN{\phi}$ is $\sqrt{d}\Lipconst^2$-Lipschitz continuous on $\R^d$ and

\item $\neuronsANN{\phi} \leq \tfrac{5}{6} 10^3 \Lipconst^4 r d^6 \epsilon^{-1}$.
\end{enumerate}
\end{proposition}

Note that the parameter estimates in \cref{ex_Lipschitz_sum_bounded,ex_Lipschitz_sum_bounded_v2} 
depend on $r$ since we approximate uniformly on the hypercube $[-r,r]^d$. However, the estimate is only 
linear in $r$, even though the volume of the hypercube is of the order $r^d$.

\begin{proposition}
\label{ex_Lipschitz_max_bounded_v3}
Let $\Lipconst,r \in [1,\infty)$ and suppose $f_i \colon \R \rightarrow \R$ is $\Lipconst$-Lipschitz 
with $|f_i(0)| \leq \Lipconst$ and $g_i \colon \R \rightarrow \R$ 1-Lipschitz with $g_i(0) = 0$, $i \in \N$.
Let $h_d \colon \R^d \rightarrow \R$, $d \in \N$, be given by $h_1 = f_1$ and 
$h_d(x) = g_d(\max\{h_{d-1}(x_1,\dots,x_{d-1}),f_d(x_d)\})$ for $d \geq 2$.
Then for all $d \in \N$ and $\epsilon \in (0,1]$, there exists 
$\phi \in \ANNs$ with $\fctReLUANN{\phi} \in C(\R^d,\R)$ such that
\begin{enumerate}[\rm 1)]
\item $\sup_{x \in [-r,r]^d} |h_d(x) - \fctReLUANN{\phi}(x)| \leq \epsilon$,

\item $\fctReLUANN{\phi}$ is $\Lipconst$-Lipschitz continuous on $\R^d$ and

\item $\neuronsANN{\phi} \leq \frac{3}{7} 10^4 \Lipconst^3 r d^{\nicefrac{13}{2}} \epsilon^{-1}$.
\end{enumerate}
\end{proposition}

Note that if the functions $g_i$, for $1 \leq i \leq d-1$, are chosen to be the identity and $g_d = f_0$, then
$h_d$ reduces to $f_0(\max\{f_1(x_1),\dots,f_d(x_d)\})$.

The functions in the previous propositions were approximated on bounded domains, but if one is willing to 
pay a higher approximation cost, one can also approximate the family of functions from, e.g., \cref{ex_Lipschitz_sum_bounded} on the entire space without the curse of dimensionality for an appropriate weight function.

\begin{proposition}
\label{ex_Lipschitz_sum_unbounded}
Let $q \in (1,\infty)$, $\Lipconst \in [1,\infty)$ and suppose 
$f_i \colon \R \rightarrow \R$, $i \in \N$, are $\Lipconst$-Lipschitz continuous with $|f_i(0)| \leq \Lipconst$.
Define $g_d \colon \R^d \rightarrow \R$ by $g_d(x) = \ssum{i=1}{d} f_i(x_i)$.
Then for all $d \in \N$ and $\epsilon \in (0,1]$, there exists $\phi \in \ANNs$ 
with $\fctReLUANN{\phi} \in C(\R^d,\R)$ such that
\begin{enumerate}[\rm 1)]
\item $\sup_{x \in \R^d} (1+\mednorm{x}^q)^{-1} |g_d(x) - \fctReLUANN{\phi}(x)| \leq \epsilon$,

\item $\fctReLUANN{\phi}$ is $\sqrt{d}\Lipconst$-Lipschitz continuous on $\R^d$ and

\item $\neuronsANN{\phi} \leq \tfrac{2}{9} 10^3 2^{3t(q + 1)} \Lipconst^{2t(q+1)} d^{t(q + 1) + 4} 
\epsilon^{-t}$ \newline for $t = \nicefrac{q}{q-1}$.
\end{enumerate}
\end{proposition}

An analogous result with the maximum function instead of the sum can be shown similarly.

In the next two propositions, we replace the sum by a product. Then we cannot
establish the approximation on an arbitrarily large domain since the Lipschitz constant 
of the product function on $[-r,r]^2$ grows linearly in $r$.

\begin{proposition}
\label{ex_Lipschitz_prod_v3}
Let $\Lipconst \in [1,\infty)$ and suppose $f_i \colon \R \rightarrow \R$ is $\Lipconst$-Lipschitz and 
$g_i \colon \R \rightarrow \R$ 1-Lipschitz with $f_i(0) = 0 = g_i(0)$, $i \in \N$.
Let $h_d \colon \R^d \rightarrow \R$, $d \in \N$, be given by $h_1 = f_1$ and 
$h_d(x) = g_d(h_{d-1}(x_1,\dots,x_{d-1}) f_d(x_d))$ for $d \geq 2$.
Set $r = \nicefrac{1}{\sqrt{8}\Lipconst}$.
Then for all $d \in \N$ and $\epsilon \in (0,\nicefrac{1}{16}]$, there exists 
$\phi \in \ANNs$ with $\fctReLUANN{\phi} \in C(\R^d,\R)$ such that
\begin{enumerate}[\rm 1)]
\item $\sup_{x \in [-r,r]^d} |h_d(x) - \fctReLUANN{\phi}(x)| \leq \epsilon$,

\item $\fctReLUANN{\phi}$ is $\Lipconst$-Lipschitz continuous on $\R^d$ and

\item $\neuronsANN{\phi} \leq \frac{3}{7} 10^4 \Lipconst^2 d^{\nicefrac{13}{2}} \epsilon^{-1}$.
\end{enumerate}
\end{proposition}

This example includes the special case $f_0(\prod_{i=1}^{d} f_i(x_i))$.
Since on large hypercubes the quantity $\CNtechLip{\CNelt_d}$, where 
$\CNelt_d$ is a catalog network representing the function $h_d$,
starts to grow exponentially in the dimension, the approximators in the proof of Proposition 
\ref{ex_Lipschitz_prod_v3} can only be built on the hypercube
$[-\nicefrac{1}{\sqrt{8}\Lipconst},\nicefrac{1}{\sqrt{8}\Lipconst}]^d$.

However, it has been shown in Proposition 3.3. of \cite{SchwabZech2019} that the product $\prod_{i=1}^{d} x_i$
can be approximated without the curse of dimensionality on the hypercube $[-1,1]^d$. 
Applying the $\log$-modification of our theory, we can recover this result and even allow for 
arbitrarily large hypercubes.

\begin{proposition}
\label{ex_prod}
Consider the functions $f_d \colon \R^d \rightarrow \R$, $d \in \N$, given 
by $f_d(x) = \prod_{i=1}^{d} x_i$, and let $r \in [1,\infty)$.
Then for all $d \in \N_{\geq 2}$ and $\epsilon \in (0,\min\{1/2,r^{-1}\}]$, there exists $\phi \in \ANNs$ with $\fctReLUANN{\phi} \in C(\R^d,\R)$ such that
\begin{enumerate}[\rm 1)]
\item $\sup_{x \in [-r,r]^d} |f_d(x) - \fctReLUANN{\phi}(x)| \leq \epsilon$,

\item $\fctReLUANN{\phi}$ is $8^{\nicefrac{(d-1)}{2}}r^{d(d-1)}$-Lipschitz continuous on $\R^d$,

\item $\neuronsANN{\phi} \leq \tfrac{1}{3} 10^5 \log_2(d) d^6 \log_2(\epsilon^{-1})$ if $r = 1$ and

\item $\neuronsANN{\phi} \leq \tfrac{2}{5} 10^5 \log_2(r) \log_2(d) d^8 \log_2(\epsilon^{-1})$ if $r \geq 2$.
\end{enumerate}
\end{proposition}

In our next result, we show the approximability of ridge functions based on a Lipschitz function.

\begin{proposition}
\label{ex_ridge}
Let $\Lipconst,r,S \in [1,\infty)$, $\theta_d \in [-S,S]^d$, $d \in \N$, and suppose 
$f \colon \R \rightarrow \R$ is $\Lipconst$-Lipschitz continuous with $|f(0)| \leq \Lipconst$.
Consider the ridge functions $g_d \colon \R^d \rightarrow \R$ given by $g_d(x) = f(\theta_d \cdot x)$.
Then for all $d \in \N$ and $\epsilon \in (0,1]$, there exists 
$\phi \in \ANNs$ with $\fctReLUANN{\phi} \in C(\R^d,\R)$ such that
\begin{enumerate}[\rm 1)]
\item $\sup_{x \in [-r,r]^d} |g_d(x) - \fctReLUANN{\phi}(x)| \leq \epsilon$,

\item $\fctReLUANN{\phi}$ is $\sqrt{d}\Lipconst S$-Lipschitz continuous on $\R^d$ and

\item $\neuronsANN{\phi} \leq \frac{1}{7} 10^3 \Lipconst^2 r S^2 d^6 \epsilon^{-1}$.
\end{enumerate}
\end{proposition}

As our last example, we consider generalized Gaussian radial basis function networks, 
i.e. weighted sums of the Gaussian function applied to the distance of $x$ to a given vector.

\begin{proposition}
\label{ex_Gaussian_RBF}
Let $N \in \N$, $r,S \in [1,\infty)$ with $r+S \geq 5$ as well as $\alpha_1,\dots,\alpha_N \in [0,S]$, $u_1,\dots,u_N \in [-S,S]$ and $v_{d,1},\dots,v_{d,N} \in [-S,S]^d$, $d \in \N$.
Consider generalized Gaussian radial basis function networks $f_d \colon \R^d \rightarrow \R$ with $N$ neurons given by $f_d(x) = \sum_{i=1}^N u_i e^{-\alpha_i \mednorm{x-v_{d,i}}^2}$.
Then for all $d \in \N$ and $\epsilon \in (0,(r+S)^{-3}]$, there exists 
$\phi \in \ANNs$ with $\fctReLUANN{\phi} \in C(\R^d,\R)$ such that
\begin{enumerate}[\rm 1)]
\item $\sup_{x \in [-r,r]^d} |f_d(x) - \fctReLUANN{\phi}(x)| \leq \epsilon$,

\item $\fctReLUANN{\phi}$ is $2(r+S)S^2 \sqrt{d} N$-Lipschitz on $\R^d$ and

\item $\neuronsANN{\phi} \leq \frac{1}{6} 10^4 (r+S) S^2 N^{\nicefrac{11}{2}} d^5 \epsilon^{-1}$.
\end{enumerate}
\end{proposition}


\section{Conclusion}
\label{section_conclusion}

It is well known that ReLU networks cannot break the curse of dimensionality for arbitrary target functions; 
see, e.g., \cite{Yarotsky2017,PeteVoigt2018,LuShenYangZhang2020}. Therefore, one has to restrict the attention to
special classes of functions. In this paper, we introduced the concept of a catalog network and proved corresponding
approximation results. Catalog networks are generalizations of standard neural networks with nonlinear activation 
functions that may change from one layer to the next as long as they belong to a given catalog. 
As such, they form a rich family of continuous functions, which offer more flexibility than other function classes 
considered in the literature. The general approximation results of \cref{section_approx_results,section_log_modification} 
give estimates on the number of parameters needed to approximate a catalog network to a given accuracy with 
ReLU-type networks. An important ingredient in our arguments is a new way to parallelize neural networks, which saves parameters 
compared to the standard parallelization. As special cases of the general results, we obtained different classes 
of functions involving sums, maxima and products of one-dimensional Lipschitz functions as well as ridge functions and 
generalized Gaussian radial basis function networks that admit a neural network approximation
without the curse of dimensionality.


\begin{appendices}

\section{Proofs for Section \ref{section_DNNs}}

\begin{IEEEproof}[\bf{Proof of \cref{lemma_parameters_of_id_conc}}]
	Abbreviate $D = \depthANN{\phi}$.
By Proposition \ref{prop_ANN_concat} and the fact that $\neuronsANN{\psi_1} = \DimANN{1}{\psi_1}(\inDimANN{\psi_1} + 1) + \DimANN{2}{\psi_1} (\DimANN{1}{\psi_1} + 1)$, we have
\begin{equation*}
	\neuronsANN{\phi \concANN \psi_1} = \neuronsANN{\phi} + \DimANN{1}{\psi_1} (\inDimANN{\phi}+1) + \DimANN{1}{\phi} (\DimANN{1}{\psi_1} - \inDimANN{\phi})
\end{equation*}
and $\DimANN{D}{\phi \concANN \psi_1} = m$.
So, by applying Proposition \ref{prop_ANN_concat} once more and observing $\DimANN{D+1}{\phi \concANN \psi_1} = \DimANN{D}{\phi} = \DimANN{2}{\psi_2}$, we obtain
\begin{equation*}
\begin{split}
	\neuronsANN{\psi_2 \concANN \phi \concANN \psi_1} &= \neuronsANN{\phi \concANN \psi_1} + \DimANN{1}{\psi_2} (\DimANN{D}{\phi}+1) + m ( \DimANN{1}{\psi_2}  - \DimANN{D}{\phi} ),
\end{split}
\end{equation*}
which completes the proof.
\end{IEEEproof}

\vspace{0.2cm}

\begin{IEEEproof}[\bf{Proof of \cref{prop_ANN_paral_NEW}}]
Assume without loss of generality that $n \geq 2$.
To simplify notation, let us introduce some abbreviations.
Write $D_j = \depthANN{\phi_j}$, $l_j^i = \DimANN{j}{\phi_i}$, $E_i = \ssum{j=1}{i} D_j$, 
$S_i = \ssum{j=1}{i} l_{D_j}^j$ and $T_i = \ssum{j=i}{n} l_0^j$.
Moreover, denote $\idconstANN_i = \idconstANN$ if $i \in \{1,\dots,n-1\}$ and $\idconstANN_i = 1$ if $i \in \{0,n\}$.
Consider the network architecture $(L_0,\dots,L_{E_n})$ of depth $E_n$ given by
\begin{equation*}
	L_k =
	\begin{cases}
		\idconstANN_i S_i + \idconstANN_i T_{i+1} &\text{if } k = E_i, \\
		\idconstANN S_{i-1} + l_m^i + \idconstANN T_{i+1} &\text{if } k = E_{i-1}+m,
	\end{cases}
\end{equation*}
where $m$ is ranging from $1$ to $D_i-1$. As discussed in 
the paragraph preceding \cref{prop_ANN_paral_NEW}, there is a skeleton $\psi \in {\cal N}$
with this architecture that realizes the parallelization of $\phi_1,\dots,\phi_n$; see also Figure \ref{fig_parallelization_novel}.
One has $\neuronsANN{\psi} = \ssum{i=1}{n} P_i$ for $P_i = \ssum{k=E_{i-1}+1}{E_i} L_k(L_{k-1}+1)$.
In the remainder of the proof, we show that $P_i \leq (\frac{11}{16}\idconstANN^2l^2n^2-1) \neuronsANN{\phi_i}$.
We distinguish the cases $D_i \geq 2$ and $D_i = 1$.
Let us begin with the former case.
By the definition of $L_k$, we have
\begin{equation*}
\begin{split}
	P_i &= (\idconstANN S_{i-1} + l_1^i + \idconstANN T_{i+1})(\idconstANN_{i-1}S_{i-1} + \idconstANN_{i-1}T_i + 1) \\
	&+ \ssum{m=2}{D_i-1} (\idconstANN S_{i-1} + l_m^i + \idconstANN T_{i+1})(\idconstANN S_{i-1} + l_{m-1}^i + \idconstANN T_{i+1} + 1) \\
	&+ (\idconstANN_iS_i + \idconstANN_iT_{i+1})(\idconstANN S_{i-1} + l_{D_i-1}^i + \idconstANN T_{i+1} + 1).
\end{split}
\end{equation*}
Now we use $\idconstANN \geq 2$, $\idconstANN \geq \idconstANN_i$, $S_i = S_{i-1} + l_{D_i}^i$, $T_i = l_0^i + T_{i+1}$ and $S_{i-1} + T_{i+1} \leq l(n-1)$, and reorder the resulting terms to obtain
\begin{equation*}
\begin{split}
	P_i &\leq \ssum{m=1}{D_i} l_m^i(l_{m-1}^i+1) + \idconstANN^2l(n-1) \ssum{m=0}{D_i} l_m^i\\
	&\quad + (\idconstANN-1) l_1^i l_0^i + (\idconstANN-1) l_{D_i}^i(l_{D_i-1}^i+1) \\
	&\quad + D_i \idconstANN l(n-1) (\idconstANN l(n-1)+1).
\end{split}
\end{equation*}
Then, we bound the second line by $2(\idconstANN-1)\neuronsANN{\phi_i}$, the sum $\ssum{m=0}{D_i} l_m^i$ by 
$\neuronsANN{\phi_i}$ and the depth $D_i$ by $\frac{1}{2}\neuronsANN{\phi_i}$ to find
\begin{equation*}
\begin{split}
	P_i &\leq \neuronsANN{\phi_i} \big[ 2\idconstANN - 1 + \idconstANN^2l(n-1) \\
	&\quad + \tfrac{1}{2} \idconstANN l(n-1)(\idconstANN l(n-1)+1) \big].
\end{split}
\end{equation*}
Finally, since $\idconstANN \geq 2$ and $n \geq 2$, the term in the brackets can be bounded by 
$\frac{11}{16}\idconstANN^2l^2n^2-1$, and the proposition follows in the case $D_i \geq 2$.
Now assume $D_i = 1$ so that $\neuronsANN{\phi_i} = l_1^i (l_0^i+1)$.
Then, by the same inequalities as in the previous case,
\begin{equation*}
\begin{split}
	P_i &= (\idconstANN_iS_i + \idconstANN_iT_{i+1})(\idconstANN_{i-1}S_{i-1} + \idconstANN_{i-1}T_i + 1) \\
	&\leq \idconstANN^2 (l(n-1) + l_1^i)(l(n-1) + l_0^i + \tfrac{1}{2}).
\end{split}
\end{equation*}
If $l_0^i = 1$, then $\neuronsANN{\phi_i} = 2 l_1^i$ and, hence,
\begin{equation*}
\begin{split}
	P_i + \neuronsANN{\phi_i} &\leq \idconstANN^2 l^2 (n-1 + \tfrac{1}{2}\neuronsANN{\phi_i})(n + \tfrac{1}{2}) + \neuronsANN{\phi_i} \\
	&\leq \idconstANN^2 l^2 \neuronsANN{\phi_i} \big[ \tfrac{n}{2}(n + \tfrac{1}{2}) + \tfrac{1}{4}\big].
\end{split}
\end{equation*}
Since $n \geq 2$, we have $\tfrac{n}{2}(n + \tfrac{1}{2}) + \tfrac{1}{4} \leq \tfrac{11}{16}n^2$, which concludes the case $l_0^i = 1$.
Finally, if $l_0^i \geq 2$, then $\neuronsANN{\phi_i} \geq 3$ and $l \geq 2$, so we obtain
\begin{equation*}
\begin{split}
	P_i + \neuronsANN{\phi_i} &\leq \idconstANN^2 l^2 ( n-1 + \tfrac{1}{2}l_1^i)(n - \tfrac{3}{4} + \tfrac{1}{2}l_0^i) + \neuronsANN{\phi_i} \\
	&\leq \idconstANN^2 l^2 \big[ (n-1)(n-\tfrac{3}{4}) + \tfrac{n-1}{2} l_0^i \\
	&\quad + (\tfrac{n-1}{2}+\tfrac{1}{8}) l_1^i + \tfrac{1}{4}l_0^il_1^i \big] + \tfrac{1}{16} \idconstANN^2 l^2 \neuronsANN{\phi_i} \\
	&\leq \idconstANN^2 l^2 \neuronsANN{\phi_i} \big[ \tfrac{1}{3}(n-1)(n-\tfrac{3}{4}) + \tfrac{n-1}{2}+\tfrac{5}{16} \big],
\end{split}
\end{equation*}
and the term in the brackets is bounded by $\tfrac{1}{3}n^2$, which finishes the last remaining case.
\end{IEEEproof}

\vspace{0.2cm}

Next, we show that the bound of \cref{prop_ANN_paral_NEW} is asymptotically sharp up to a constant by computing the number of parameters of an $n$-fold diagonal parallelization of the same skeleton $\phi_1 = \dots = \phi_n$ satisfying $\depthANN{\phi_1} \geq 2$ and having a single neuron in all its layers in the case $\idconstANN = 2$.
With the notation from the previous proof, we have $S_i = i$ and $T_i = n-i+1$.
Thus, the formula for $P_i$ reads
\begin{equation*}
\begin{split}
	P_i &= (2n-1)2n(D_i-2) \\
	&\quad +
	\begin{cases}
		(2n-1)(n+1) + 4n^2, &\text{if } i = 1, \\
		(2n-1)(2n+1) + 4n^2, &\text{if } 2 \leq i \leq n-1, \\
		(2n-1)(2n+1) + 2n^2, &\text{if } i = n.
	\end{cases}
\end{split}
\end{equation*}
Since $\neuronsANN{\phi_1} = 2D_i$, we find for the diagonal parallelization
\begin{equation*}
	\neuronsANN{\paralANN{\IDANN{}}{\phi_1,\dots,\phi_n}} = P_1 + (n-2)P_2 + P_n =	(2n^3-n^2)\neuronsANN{\phi_1}.
\end{equation*}
Together with the upper bound, this proves the asymptotic sharpness of \cref{prop_ANN_paral_NEW} up to a factor of at most $\tfrac{11}{8}$.

\vspace{0.2cm}

\begin{IEEEproof}[\bf{Proof of \cref{cor_parameters_of_id_conc}}]
	Abbreviate $D = \depthANN{\phi}$, $k = \inDimANN{\phi}$ and $n = \outDimANN{\phi}$. 
	First, assume $D \geq 2$.
Lemma \ref{lemma_parameters_of_id_conc} yields
\begin{equation*}
\begin{split}
	\neuronsANN{\IDANN{n} \concANN \phi \concANN \IDANN{k}} &= \neuronsANN{\phi} + \DimANN{1}{\IDANN{k}} (k+1) + \DimANN{1}{\IDANN{n}} (n+1) \\
	&\quad + \DimANN{1}{\phi} (\DimANN{1}{\IDANN{k}} - k) + \DimANN{D-1}{\phi}(\DimANN{1}{\IDANN{n}} - n).
\end{split}
\end{equation*}
Note that $\DimANN{1}{\IDANN{k}}$ and $\DimANN{1}{\IDANN{n}}$ are at most $\idconstANN k$ and $\idconstANN n$, respectively.
This and the fact that $\DimANN{1}{\phi} + \DimANN{D-1}{\phi} \leq \frac{2}{3}\neuronsANN{\phi}$ imply
\begin{equation*}
\begin{split}
	\neuronsANN{\IDANN{n} \concANN \phi \concANN \IDANN{k}} &\leq \neuronsANN{\phi} + 2\idconstANN m(m + 1) \\
	&\quad + (\DimANN{1}{\phi}+\DimANN{D-1}{\phi}) (\idconstANN m - 1) \\
	&\leq \tfrac{5}{6} \idconstANN m \neuronsANN{\phi} + 2 \idconstANN^2 m^2,
\end{split}
\end{equation*}
where the last inequality holds because $\idconstANN \geq 2$.
Now, suppose $D = 1$.
Then, by Lemma \ref{lemma_parameters_of_id_conc},
\begin{equation*}
\begin{split}
	\neuronsANN{\IDANN{n} \concANN \phi \concANN \IDANN{k}} &\leq n + \idconstANN k(k+1) + \idconstANN n(n+1) + \idconstANN^2 n k \\
	&\leq \tfrac{5}{6} \idconstANN \neuronsANN{\phi} + \idconstANN^2 m^2 + \tfrac{7}{6}\idconstANN m (m+1) + m \\
	&\leq \tfrac{5}{6} \idconstANN m \neuronsANN{\phi} + \tfrac{29}{12} \idconstANN^2 m^2,
\end{split}
\end{equation*}
where we again used $\idconstANN \geq 2$.
\end{IEEEproof}


\section{Proofs for Section \ref{section_CNs}}

\begin{IEEEproof}[\bf{Proof of \cref{ex_weights_OoG}}]
Denote $s = \max\{b_0,\dots,b_q\}$. Since the coefficients $a_0,\dots,a_q$ of $g$ 
are non-negative, one has $g(rx) \leq r^s g(x)$ for all $x \in [0,\infty)$ and $r \in [1,\infty)$.
This and the assumption that $f$ is non-decreasing yield 
\begin{equation*}
\weight(x) \leq \frac{f(rx)}{\max\{g(x),\delta\}} \leq \frac{f(rx) r^s}{\max\{g(rx),\delta\}} = r^s \weight(rx)
\end{equation*}
for all $x \in [0,\infty)$, $r \in [1,\infty)$. That $f$ is non-decreasing also gives
\begin{equation*}
	\weight(x) \leq \frac{f(1)}{\max\{g(x),\delta\}} \leq \frac{\max\{g(1),\delta\}}{\delta} \weight(1)
\end{equation*}
for all $x \in [0,1)$. Combining the previous two estimates yields
\begin{equation*}
	\frac{\delta}{\max\{g(1),\delta\}} \weight(x) \leq \weight(\max\{x,1\}) \leq r^s \weight(r\max\{x,1\})
\end{equation*}
for all $x \in [0,\infty)$, $r \in [1,\infty)$, which finishes the proof.
\end{IEEEproof}

\vspace{0.2cm}

\begin{IEEEproof}[\bf{Proof of \cref{lemma_main_layer_Lipschitz}}]
Assume $\CNelt$ is of the form $\CNelementfixed$. As discussed after
\cref{def_approx_catalog}, every $f \in \Catalog$ is $\Lip_f$-Lipschitz continuous on the set $\approxSet{f}$.
For $k \in \{1,\dots,D\}$, $j \in \{1,\dots,n_k\}$ and $x \in \R^{l_{2k-1}}$, denote by $x_{(k,j)}$ 
the vector $( x_{\CatDimD{f_{k,1}} + \cdots + \CatDimD{f_{k,j-1}}+1} , \dots, x_{\CatDimD{f_{k,1}} + \cdots + \CatDimD{f_{k,j}}} ) \in \R^{\CatDimD{f_{k,j}}}$.
Then, for all $k \in \{1,\dots,D\}$ and $x,y \in \CNrdlayer{k}{\CNelt}$,
\begin{equation}
\label{norm_remark_lem_0}
\begin{split}
	&\mednorm{\CNnonlin{k}{\CNelt}(x) - \CNnonlin{k}{\CNelt}(y)}^2 \\
	&= \ssum{j=1}{n_k} \mednorm{f_{k,j}(x_{(k,j)}) - f_{k,j}(y_{(k,j)})}^2 \\
	&\leq \ssum{j=1}{n_k} \Lip_{f_{k,j}}^2 \mednorm{x_{(k,j)} - y_{(k,j)}}^2 \leq [\CNlayerLip{k}{\CNelt}]^2 \mednorm{x-y}^2,
\end{split}
\end{equation}
which is what we wanted to show.
\end{IEEEproof}


\section{Proofs for Section \ref{section_ex_approx_cat}}

\begin{IEEEproof}[\bf{Proof of \cref{lemma_sq_mirror}}]
Choose $\phi_{r,1}, \phi_{r,2} \in \ANNs$ of depth 1 such that $\fctReLUANN{\phi_{r,1}} \in C(\R,\R^2)$ realizes 
$x \mapsto (\frac{x}{r}, -\frac{x}{r})$ and $\fctReLUANN{\phi_{r,2}} \in C(\R^2,\R)$ realizes $(x,y) \mapsto r^2(x+y)$.
If $(\phi_{\epsilon})_{\epsilon \in (0,1]} \subseteq \ANNs$ are the $\varepsilon$-approximations of the
square function on $[0,1]$ derived in Proposition 3.3 of \cite{GrohsHorJenZimm2019},
then\footnote{Here, we understand $\phi_{r^{-2}\epsilon}$ as $\phi_1$ if $r^{-2}\epsilon > 1$.} the ReLU-realization of 
$\psi_{r,\epsilon} = \phi_{r,2} \concANN \paralANN{}{\phi_{r^{-2}\epsilon},\phi_{r^{-2}\epsilon}} \concANN \phi_{r,1}$ approximates the square function on $[-r,r]$ with accuracy $\epsilon$.
To see this, note that $\fctReLUANN{\psi_{r,\epsilon}}(x) = r^2 \fctReLUANN{\phi_{r^{-2}\epsilon}}(\tfrac{|x|}{r})$ for all $x \in \R$ since $\phi_{r^{-2}\epsilon} = \ReLU$ on $\R \backslash [0,1]$.
This also implies $\fctReLUANN{\psi_{r,\epsilon}}(x) = r|x|$ for all $x \in \R \backslash [-r,r]$ as well as the 
$2r$-Lipschitz continuity. Finally, (\ref{prop_ANN_concat_item_depth_1_1}) and (\ref{prop_ANN_concat_item_depth_1_2}) of \cref{prop_ANN_concat} together with (\ref{prop_ANN_paral_same_depth_item_same_network}) of \cref{prop_ANN_paral_same_depth} 
assure that $\neuronsANN{\psi_{r,\epsilon}} \leq 4 \neuronsANN{\phi_{r^{-2}\epsilon}}$, which concludes the proof.
\end{IEEEproof}

\vspace{0.2cm}

\begin{IEEEproof}[\bf{Proof of \cref{prop_product_approx}}]
Choose $\psi_1, \psi_2 \in \ANNs$ of depth 1 such that $\fctReLUANN{\psi_1} \in C(\R^2,\R^2)$ realizes $(x,y) \mapsto (x+y,x-y)$ and  
$\fctReLUANN{\psi_2} \in C(\R^2,\R)$ realizes $(x,y) \mapsto \tfrac{1}{4}(x-y)$.
If $(\psi_{r,\epsilon})_{r \in (0,\infty),\epsilon \in (0,1]} \subseteq \ANNs$ denote the $\varepsilon$-approximations of the square 
function on the interval $[-r,r]$ from \cref{lemma_sq_mirror}, then the ReLU-realization of $\chi_{r,\epsilon} = \psi_2 \concANN \paralANN{}{\psi_{2r,2\epsilon},\psi_{2r,2\epsilon}} \concANN \psi_1$ approximates the product function on $[-r,r]^2$ with accuracy $\epsilon$.
Furthermore, $\fctReLUANN{\chi_{r,\epsilon}}$ is $\sqrt{8}r$-Lipschitz continuous because $\fctReLUANN{\psi_{2r,2\epsilon}}$ is $4r$-Lipschitz continuous and, hence,
\begin{equation*}
\begin{split}
	&|\fctReLUANN{\chi_{r,\epsilon}}(x_1,x_2) - \fctReLUANN{\chi_{r,\epsilon}}(y_1,y_2)|  \\
	&\leq 2r ( |x_1-y_1| + |x_2-y_2| ) \leq \sqrt{8}r \mednorm{(x_1-y_1,x_2-y_2)}.
\end{split}
\end{equation*}
As in the proof of \cref{lemma_sq_mirror}, combining (\ref{prop_ANN_concat_item_depth_1_1}) and (\ref{prop_ANN_concat_item_depth_1_2}) of \cref{prop_ANN_concat} with (\ref{prop_ANN_paral_same_depth_item_same_network}) of \cref{prop_ANN_paral_same_depth},
shows that $\neuronsANN{\chi_{r,\epsilon}} \leq 4 \neuronsANN{\psi_{2r,2\epsilon}}$, from which the proposition follows.
\end{IEEEproof}


\section{Proofs for Section \ref{section_approx_results}}

\begin{IEEEproof}[\bf{Proof of \cref{lemma_main_layer_approx}}]
Assume $\CNelt$ is of the form $\CNelementfixed$. Fix any $k \in \{1,\dots,D\}$ and $\delta \in (0,\epsilon]$.
The assumption that $\Catalog$ is $\ApproxCat{\epsilon}{\kappa}$-approximable for $\kappa = (\kappa_0,\kappa_1,\kappa_2,\kappa_3)$
guarantees that there exist skeletons $\psi_{j} \in \ANNs$, $j \in \{1,\dots,n_k\}$, 
such that the $a$-realizations $\fctANN{\psi_{j}} \in C(\R^{\CatDimD{f_{k,j}}},\R^{\CatDimT{f_{k,j}}})$ satisfy
\begin{enumerate}[\rm 1)]
\item $\mednorm{\fctANN{\psi_{j}}(x) - \fctANN{\psi_{j}}(y)} \leq \Lip_{f_{k,j}} \mednorm{x-y}$ for all $x,y \in \R^{\CatDimD{f_{k,j}}}$,

\item $\weight(\mednorm{x}) \mednorm{f_{k,j}(x) - \fctANN{\psi_{j}}(x)} \leq \frac{\delta}{\sqrt{n_k}}$ for all $x \in \approxSet{f_{k,j}}$

\item and $\neuronsANN{\psi_{j}} \leq \kappa_1 \max\{\CatDimD{f_{k,j}},\CatDimT{f_{k,j}}\}^{\kappa_2} n_k^{\nicefrac{\kappa_3}{2}} \delta^{-\kappa_3}$.
\end{enumerate}
Pick an $\IDANN{} \in \ANNs$ such that $a$ fulfills the $\idconstANN$-identity requirement with $\IDANN{}$, 
and let $\phi \in \ANNs$ be the $I$-parallelization $\phi = \paralANN{\IDANN{}}{\psi_{1},\dots,\psi_{n_k}}$.
For $j \in \{1,\dots,n_k\}$ and $x \in \R^{l_{2k-1}}$, denote $x_{(k,j)} = ( x_{\CatDimD{f_{k,1}} + \cdots + \CatDimD{f_{k,j-1}}+1} , \dots, x_{\CatDimD{f_{k,1}} + \cdots + \CatDimD{f_{k,j}}} )$.
Then, for all $x,y \in \R^{l_{2k-1}}$, 
\begin{equation}
\label{norm_remark_lem_1}
\begin{split}
	&\mednorm{\fctANN{\phi}(x) - \fctANN{\phi}(y)}^2 \\
	&= \ssum{j=1}{n_k} \mednorm{\fctANN{\psi_{j}}(x_{(k,j)}) - \fctANN{\psi_{j}}(y_{(k,j)})}^2 \\
	&\leq \ssum{j=1}{n_k} \Lip_{f_{k,j}}^2 \mednorm{x_{(k,j)} - y_{(k,j)}}^2 \leq [\CNlayerLip{k}{\CNelt}]^2 \mednorm{x-y}^2.
\end{split}
\end{equation}
Moreover, since $\weight$ is non-increasing, we obtain, for all $x \in \CNrdlayer{k}{\CNelt}$,
\begin{equation}
\label{norm_remark_lem_2}
\begin{split}
	\mednorm{\CNnonlin{k}{\CNelt}(x) &- \fctANN{\phi}(x)}^2 \\
	&= \ssum{j=1}{n_k} \mednorm{f_{k,j}(x_{(k,j)}) - \fctANN{\psi_{j}}(x_{(k,j)})}^2 \\
	&\leq \frac{\delta^2}{n_k} \ssum{j=1}{n_k} [\weight( \mednorm{x_{(k,j)}} )]^{-2} \leq \delta^2 [\weight( \mednorm{x} )]^{-2}.
\end{split}
\end{equation}
It remains to estimate the number of parameters $\neuronsANN{\phi}$.
Since $\inDimANN{\psi_j} = \CatDimD{f_{k,j}} \leq l_{2k-1}$ and $\outDimANN{\psi_j} = \CatDimT{f_{k,j}} \leq l_{2k}$ for each $j \in \{1,\dots,n_k\}$, \cref{prop_ANN_paral_NEW} yields
\begin{equation}
\label{exponent_remark_lem}
\begin{split}
	\neuronsANN{\phi} &\leq \tfrac{11}{16} \idconstANN^2 \max\{l_{2k-1},l_{2k}\}^2 n_k^2 \ssum{j=1}{n_k} \neuronsANN{\psi_j} \\
	&\leq \tfrac{11}{16} \kappa_1 \idconstANN^2 \max\{l_{2k-1},l_{2k}\}^{\kappa_2 + 2} n_k^{3+\frac{\kappa_3}{2}} \delta^{-\kappa_3}.
\end{split}
\end{equation}
Note that we always have $n_k \leq \max\{l_{2k-1},l_{2k}\}$, which yields (\ref{lemma_main_layer_approx_item}).
If $\CatDimD{f} \leq d$ and $\CatDimT{f} \leq d$ for some $d \in \N$ and all $f \in \Catalog$, then we use the estimate $\inDimANN{\psi_j} = \CatDimD{f_{k,j}} \leq d$ instead of $\inDimANN{\psi_j} = \CatDimD{f_{k,j}} \leq l_{2k-1}$ (and similarly for $\outDimANN{\psi_j}$) to 
obtain $\neuronsANN{\phi} \leq \frac{11}{16} \idconstANN^2 d^2 n_k^2 \ssum{j=1}{n_k} \neuronsANN{\psi_j}$ in (\ref{exponent_remark_lem}), which shows (\ref{lemma_main_layer_approx_item_prime}).
\end{IEEEproof}

\vspace{0.2cm}

\begin{IEEEproof}[\bf{Proof of Theorem \ref{theorem_main}}]
We split the proof into two parts. In the first part, we construct an approximating neural network and bound the approximation error.
In the second part, we estimate the number of parameters of the network.
Assume $\CNelt \in \CNfixed$ is of the form $\CNelt = \CNelementfixed$. Denote
$\auxproofcomplayerreal{0} = \ID{l_0}$ and
\begin{equation*}
	\auxproofcomplayerreal{k} = \CNnonlin{k}{\CNelt} \circ \CNaff{k}{\CNelt} \circ \CNnonlin{k-1}{\CNelt} \circ \cdots \circ \CNnonlin{1}{\CNelt} \circ \CNaff{1}{\CNelt}
\end{equation*}
for $k \in \{1,\dots,D\}$. Before constructing the approximating network, we show by induction over $k$ that 
\begin{equation}
\label{proof_thrm_induc_1}
\begin{split}
	&\mednorm{(\CNaff{k}{\CNelt} \circ \auxproofcomplayerreal{k-1})(x)} \leq \opnorm{V_k} \Big( \sprod{j=1}{k-1} \CNlayerLip{j}{\CNelt} \opnorm{V_j} \Big) \mednorm{x} + \mednorm{b_k} \\
	&\quad + \ssum{i=1}{k-1} \opnorm{V_k} \Big( \sprod{j=i+1}{k-1} \CNlayerLip{j}{\CNelt} \opnorm{V_j} \Big) \big( \CNlayerLip{i}{\CNelt} \mednorm{b_i} + \mednorm{\CNnonlin{i}{\CNelt}(0)} \big)
\end{split}
\end{equation}
for all  $k \in \{1,\dots,D\}$ and $x \in \CNredom{\CNelt}$. The base case $k=1$ reduces to 
the obvious inequality $\mednorm{\CNaff{1}{\CNelt}(x)} \leq \opnorm{V_1} \mednorm{x} + \mednorm{b_1}$.
For the induction step, suppose the claim is true for a given $k \in \{1,\dots,D-1\}$.
Then we obtain from \cref{lemma_main_layer_Lipschitz} that
\begin{equation*}
\begin{split}
	\mednorm{(\CNaff{k+1}{\CNelt} &\circ \auxproofcomplayerreal{k})(x)} \\
	&\leq \opnorm{V_{k+1}} \mednorm{(\CNnonlin{k}{\CNelt} \circ \CNaff{k}{\CNelt} \circ \auxproofcomplayerreal{k-1})(x)} + \mednorm{b_{k+1}} \\
	&\leq \opnorm{V_{k+1}} \CNlayerLip{k}{\CNelt} \mednorm{(\CNaff{k}{\CNelt} \circ \auxproofcomplayerreal{k-1})(x)} \\
	&\quad + \opnorm{V_{k+1}} \mednorm{\CNnonlin{k}{\CNelt}(0)} + \mednorm{b_{k+1}}.
\end{split}
\end{equation*}
for all $x \in \CNredom{\CNelt}$, where we used that the sets $\approxSet{f}$ contain $0$.
Applying the induction hypothesis to $\mednorm{(\CNaff{k}{\CNelt} \circ \auxproofcomplayerreal{k-1})(x)}$, 
we obtain that (\ref{proof_thrm_induc_1}) holds for $k+1$. Next, observe that 
\begin{equation}
\label{norm_remark_thrm}
	\mednorm{\CNnonlin{k}{\CNelt}(0)}^2 = \ssum{j=1}{n_k} \mednorm{f_{k,j}(0)}^2 \leq \kappa_0^2 n_k \leq \kappa_0^2 \CNwidth{\CNelt}
\end{equation}
for all $k \in \{1,\dots,D\}$. Hence, (\ref{proof_thrm_induc_1}) yields 
\begin{equation}
\label{proof_thrm_precision_weight_input}
\begin{split}
	\mednorm{(\CNaff{D}{\CNelt} &\circ \auxproofcomplayerreal{D-1})(x)} \\
	&\leq \CNtechLip{\CNelt} \mednorm{x} + \CNtransl{\CNelt} + D \CNtechLip{\CNelt} \big( \CNtransl{\CNelt} + \kappa_0\sqrt{\CNwidth{\CNelt}} \big) \\
	&\leq 4 \kappa_0 \CNtransl{\CNelt} \CNdepth{\CNelt} \sqrt{\CNwidth{\CNelt}} \CNtechLip{\CNelt} \max\{1,\mednorm{x}\}
\end{split}
\end{equation}
for all $x \in \CNredom{\CNelt}$. 

To construct an approximating network, we note that it follows from
Lemma \ref{lemma_main_layer_approx} that for all $\delta \in (0,\epsilon]$ and $k \in \{1,\dots,D\}$, 
there exists $\psi_{\delta,k} \in \ANNs$ such that the $a$-realization $\fctANN{\psi_{\delta,k}} \in C(\R^{l_{2k-1}}\R^{l_{2k}})$ satisfies
\begin{enumerate}[\rm 1)]

\item $\mednorm{\fctANN{\psi_{\delta,k}}(x) - \fctANN{\psi_{\delta,k}}(y)} \leq \CNlayerLip{k}{\CNelt} \mednorm{x-y}$ for all $x,y \in \R^{l_{2k-1}}$,

\item $\weight(\mednorm{x}) \mednorm{\CNnonlin{k}{\CNelt}(x) - \fctANN{\psi_{\delta,k}}(x)} \leq \delta$ for all $x \in \CNrdlayer{k}{\CNelt}$ and

\item\label{proof_thrm_main_item_neurons} $\neuronsANN{\psi_{\delta,k}} \leq \frac{11}{16} \kappa_1 \idconstANN^2 \CNwidth{\CNelt}^{\kappa_2 + \frac{\kappa_3}{2} + 5} \delta^{-\kappa_3}$.
\end{enumerate}
Moreover, since each $\CNaff{k}{\CNelt}$ is an affine function, there exist unique $\chi_k \in \ANNs$ of depth 1, such that 
$\fctANN{\chi_k} = \CNaff{k}{\CNelt}$ for all $k \in \{1,\dots,D\}$.
Let $\varphi_{\delta,k} \in \ANNs$ be given by $\varphi_{\delta,k} = \psi_{\delta,k} \concANN \chi_k \concANN \cdots \concANN \psi_{\delta,1} \concANN \chi_1$.
The $a$-realization of $\varphi_{\delta, D}$ will be our approximation network.
To verify that it does the job in terms of the approximation precision, we show that 
\begin{equation}
\label{proof_thrm_induc_2}
	\mednorm{\auxproofcomplayerreal{k}(x) - \fctANN{\varphi_{\delta,k}}(x)} \leq \sum_{i=1}^{k} \frac{\delta \sprod{j=i+1}{k} \CNlayerLip{j}{\CNelt} \opnorm{V_j}}{\weight(\mednorm{\CNaff{i}{\CNelt} \circ \auxproofcomplayerreal{i-1}(x)})}
\end{equation}
for all $k \in \{1,\dots,D\}$, $\delta \in (0,\epsilon]$ and $x \in \CNredom{\CNelt}$ by induction over $k$.
The base case $k=1$ holds by the approximation property of $\psi_{\delta,1}$ and the fact that $\CNaff{1}{\CNelt}(x) \in \CNrdlayer{1}{\CNelt}$ for all $x \in \CNredom{\CNelt}$.
For the induction step, we assume \eqref{proof_thrm_induc_2} holds for a given $k \in \{1,\dots,D-1\}$.
By the Lipschitz and approximation properties of $\psi_{\delta,k+1}$, we obtain 
\begin{equation*}
\begin{split}
	\mednorm{\auxproofcomplayerreal{k+1}(x) &- \fctANN{\varphi_{\delta,k+1}}(x)} \leq \delta [\weight(\mednorm{(\CNaff{k+1}{\CNelt} \circ \auxproofcomplayerreal{k})(x)}) ]^{-1} \\
	&\qquad + \CNlayerLip{k+1}{\CNelt} \opnorm{V_{k+1}} \mednorm{\auxproofcomplayerreal{k}(x) - \fctANN{\varphi_{\delta,k}}(x)}
\end{split}
\end{equation*}
for all $\delta \in (0,\epsilon]$ and $x \in \CNredom{\CNelt}$,
where we used that $(\CNaff{k+1}{\CNelt} \circ \auxproofcomplayerreal{k})(x) \in \CNrdlayer{k+1}{\CNelt}$ for all $x \in \CNredom{\CNelt}$.
Using the induction hypothesis on $\mednorm{\auxproofcomplayerreal{k}(x) - \fctANN{\varphi_{\delta,k}}(x)}$, we obtain 
that (\ref{proof_thrm_induc_2}) holds for $k+1$.
Now, we combine (\ref{proof_thrm_precision_weight_input}), (\ref{proof_thrm_induc_2}) and the 
assumption that the decay of $\weight$ is controlled by $(\OoGone,\OoGtwo)$ to find that
\begin{equation}
\label{proof_thrm_precision_final}
\begin{split}
	&\mednorm{\auxproofcomplayerreal{D}(x) - \fctANN{\varphi_{\delta,D}}(x)} \\
	&\leq \delta D \CNtechLip{\CNelt} \big[\weight\big(4 \kappa_0 \CNtransl{\CNelt} \CNdepth{\CNelt} \sqrt{\CNwidth{\CNelt}} \CNtechLip{\CNelt} \max\{1,\mednorm{x}\}\big) \big]^{-1} \\
	&\leq \OoGone \delta [4 \kappa_0 \CNtransl{\CNelt}]^{\OoGtwo} [\CNdepth{\CNelt} \CNtechLip{\CNelt}]^{\OoGtwo+1} \CNwidth{\CNelt}^{\OoGtwo/2} [\weight(\mednorm{x})]^{-1}
\end{split}
\end{equation}
for all $\delta \in (0,\epsilon]$ and $x \in \CNredom{\CNelt}$.

Next note that it follows from the fact that the concatenation of Lipschitz functions is again Lipschitz with constant 
equal to the product of the original Lipschitz constants that $\fctANN{\varphi_{\delta,D}}$ is $\CNfullLip{\CNelt}$-Lipschitz.

It remains to estimate the number of parameters of the constructed network. To do this, we slightly modify 
$\varphi_{\delta, D}$ by interposing identity networks. This does not change the realization but reduces the 
worst-case parameter count, as discussed before \cref{cor_parameters_of_id_conc}.
Choose $\IDANN{} \in \ANNs$ for which $a$ fulfills the $\idconstANN$-identity requirement and let $\IDANN{d} = \paralANN{}{\IDANN{},\dots,\IDANN{}}$, $d \in \N$. Define $\rho_{\delta,k} = \IDANN{l_{2k}} \concANN \psi_{\delta,k} \concANN \IDANN{l_{2k-1}}$ for $k \in \{1,\dots,D\}$ and $\delta \in (0,\epsilon]$. 
Combining \cref{cor_parameters_of_id_conc} with our parameter bound for $\neuronsANN{\psi_{\delta,k}}$, 
we obtain 
\begin{equation}
\label{proof_thrm_cost_layer}
\begin{split}
	\neuronsANN{\rho_{\delta,k}} &\leq \tfrac{55}{96} \kappa_1 \idconstANN^3 \CNwidth{\CNelt}^{\kappa_2 + \kappa_3/2 + 6} \delta^{-\kappa_3} + \tfrac{29}{12} \idconstANN^2 \CNwidth{\CNelt}^2 \\
	&\leq \tfrac{57}{32} \kappa_1 \idconstANN^3 \CNwidth{\CNelt}^{\kappa_2 + \kappa_3/2 + 6} \delta^{-\kappa_3}
\end{split}
\end{equation}
for all $\delta \in (0,\epsilon]$ and $k \in \{1,\dots,D\}$, where we used that $\idconstANN \geq 2$.
Since $\DimANN{1}{\IDANN{d}} \leq \idconstANN d$ for all $d \in \N$,
\cref{prop_ANN_concat} yields that 
\begin{equation}
\label{proof_thrm_cost_ext_layer}
	\neuronsANN{\rho_{\delta,k} \concANN \chi_k} \leq \neuronsANN{\rho_{\delta,k}} + \idconstANN \CNwidth{\CNelt}^2
\end{equation}
for all $\delta \in (0,\epsilon]$ and $k \in \{1,\dots,D\}$. Next, we show that 
\begin{equation}
\begin{split}
\label{proof_thrm_cost_induc}
	\neuronsANN{\rho_{\delta,k} &\concANN \chi_k \concANN \cdots \concANN \rho_{\delta,1} \concANN \chi_1} \\
	&\leq  (k-1)\idconstANN^2 \CNwidth{\CNelt}^2 + \ssum{j=1}{k} \neuronsANN{\rho_{\delta,j} \concANN \chi_j}
\end{split}
\end{equation}
for all $k \in \{1,\dots,D\}$ and $\delta \in (0,\epsilon]$ by induction over $k$.
The base case $k=1$ is trivially satisfied. For the induction step, we assume \eqref{proof_thrm_cost_induc} 
holds for $k \in \{1,\dots,D-1\}$. Then \cref{prop_ANN_concat} and the induction hypothesis show that 
\begin{equation*}
\begin{split}
	\neuronsANN{&\rho_{\delta,k+1} \concANN \chi_{k+1} \concANN \cdots \concANN \rho_{\delta,1} \concANN \chi_1} \\
	&\leq \neuronsANN{\rho_{\delta,k+1} \concANN \chi_{k+1}} + \neuronsANN{\rho_{\delta,k} \concANN \chi_k \concANN \cdots \concANN \rho_{\delta,1} \concANN \chi_1} \\
	&\quad + \idconstANN^2 l_{2(k+1)-1} l_{2k} \\
	&\leq k\idconstANN^2 \CNwidth{\CNelt}^2 + \ssum{j=1}{k+1} \neuronsANN{\rho_{\delta,j} \concANN \chi_j}
\end{split}
\end{equation*}
for all $\delta \in (0,\epsilon]$, which completes the induction.
Now we combine (\ref{proof_thrm_cost_ext_layer}) and (\ref{proof_thrm_cost_induc}) to obtain 
\begin{equation*}
\begin{split}
	\neuronsANN{\rho_{\delta,D} \concANN \chi_D \concANN \cdots &\concANN \rho_{\delta,1} \concANN \chi_1} \leq \tfrac{3}{2} \idconstANN^2 D \CNwidth{\CNelt}^2 + \ssum{j=1}{D} \neuronsANN{\rho_{\delta,j}}
\end{split}
\end{equation*}
for all $\delta \in (0,\epsilon]$, where we used $\idconstANN \geq 2$ again.
Using (\ref{proof_thrm_cost_layer}) gives
\begin{equation}
\label{proof_thrm_cost_final}
\begin{split}
	\neuronsANN{\rho_{\delta,D} &\concANN \chi_D \concANN \cdots \concANN \rho_{\delta,1} \concANN \chi_1} \\
	&\leq \tfrac{3}{2} \idconstANN^2 D \CNwidth{\CNelt}^2  + \tfrac{57}{32} D \kappa_1 \idconstANN^3 \CNwidth{\CNelt}^{\kappa_2 + \kappa_3/2 + 6} \delta^{-\kappa_3} \\
	&\leq \tfrac{81}{32} \kappa_1 \idconstANN^3 \CNdepth{\CNelt} \CNwidth{\CNelt}^{\kappa_2 + \kappa_3/2 + 6} \delta^{-\kappa_3}
\end{split}
\end{equation}
for all $\delta \in (0,\epsilon]$.

We conclude by summarizing what we have proved so far.
Motivated by (\ref{proof_thrm_precision_final}), let $\eta \in (0,\epsilon]$ be given by
\begin{equation*}
	\eta = \epsilon \big( \OoGone [4 \kappa_0 \CNtransl{\CNelt}]^{\OoGtwo} [\CNdepth{\CNelt} \CNtechLip{\CNelt}]^{\OoGtwo+1} \CNwidth{\CNelt}^{\OoGtwo/2} \big)^{-1},
\end{equation*}
and define $\phi = \rho_{\eta,D} \concANN \chi_D \concANN \cdots \concANN \rho_{\eta,1} \concANN \chi_1$. Then, $\fctANN{\phi} = \fctANN{\varphi_{\eta,D}}$ and $\auxproofcomplayerreal{D} = \CNfct{\CNelt}$.
So, one obtains from (\ref{proof_thrm_precision_final}) that
\begin{equation*}
	\sup\nolimits_{x \in \CNredom{\CNelt}} \weight(\mednorm{x}) \mednorm{\CNfct{\CNelt}(x) - \fctANN{\phi}(x)} \leq \epsilon,
\end{equation*}
and (\ref{proof_thrm_cost_final}) gives
\begin{equation*}
\begin{split}
	\neuronsANN{\phi} &\leq \tfrac{81}{32} \idconstANN^3 \kappa_1 \OoGone^{\kappa_3} [4 \kappa_0 \CNtransl{\CNelt}]^{\kappa_3\OoGtwo} \CNtechLip{\CNelt}^{\kappa_3(\OoGtwo + 1)} \\
	&\quad \times \CNdepth{\CNelt}^{\kappa_3(\OoGtwo+1)+1} \CNwidth{\CNelt}^{\kappa_2 + \frac{\kappa_3}{2}(\OoGtwo+1) + 6} \epsilon^{-\kappa_3},
\end{split}
\end{equation*}
which completes the proof.
\end{IEEEproof}

\vspace{0.2cm}

\begin{IEEEproof}[\bf{Proof of \cref{rem_norm}}]
To see that the proof of \cref{theorem_main} also works for the sup-norm, we note that 
in (\ref{norm_remark_lem_0}), (\ref{norm_remark_lem_1}), (\ref{norm_remark_lem_2}) and (\ref{norm_remark_thrm}),
we used the property $\mednorm{x}^2 = \ssum{j=1}{n_k} \mednorm{x_{(k,j)}}^2$ of the Euclidean norm.
However, since the sup-norm satisfies $\mednorm{x}^2_{\infty} = \max_{j \in \{1,\dots,n_k\}} \mednorm{x_{(k,j)}}^2_{\infty}$, 
and we did not use any specific property of the Euclidean norm anywhere else in the proof, 
\cref{theorem_main} still holds for the sup-norm.
\end{IEEEproof}


\section{Proofs for Section \ref{section_log_modification}}

\begin{IEEEproof}[\bf{Proof of \cref{theorem_main_LOG}}]
In case the catalog ${\cal F}$ in \cref{lemma_main_layer_approx} is $\ApproxCat{\epsilon}{\kappa}$-log-approximable
instead of $\ApproxCat{\epsilon}{\kappa}$-approximable, a slight modification of the proof yields a version 
of \cref{lemma_main_layer_approx}, for which the parameter bound in (\ref{lemma_main_layer_approx_item}) is
\begin{equation*}
	\tfrac{11}{16} \kappa_1 \idconstANN^2 \max\{l_{2k-1},l_{2k}\}^{\kappa_2 + 5} \Big[\!\log_2\!\Big(\frac{\sqrt{\max\{l_{2k-1},l_{2k}\}}}{\delta}\Big)\Big]^{\kappa_3}.
\end{equation*}
Then, in the proof of \cref{theorem_main}, one only needs to replace (\ref{proof_thrm_main_item_neurons}) with 
\begin{equation*}
	\neuronsANN{\psi_{\delta,k}} \leq \tfrac{11}{16} \kappa_1 \idconstANN^2 \CNwidth{\CNelt}^{\kappa_2 + 5} \big[\!\log_2\!\big( \sqrt{\CNwidth{\CNelt}} \delta^{-1} \big)\big]^{\kappa_3},
\end{equation*}
(\ref{proof_thrm_cost_layer}) with
\begin{equation*}
	\neuronsANN{\rho_{\delta,k}} \leq \tfrac{57}{32} \kappa_1 \idconstANN^3 \CNwidth{\CNelt}^{\kappa_2 + 6} \big[\!\log_2\!\big( \sqrt{\CNwidth{\CNelt}} \delta^{-1} \big)\big]^{\kappa_3}
\end{equation*}
and (\ref{proof_thrm_cost_final}) with
\begin{equation*}
\begin{split}
	\neuronsANN{\rho_{\delta,D} &\concANN \chi_D \concANN \cdots \concANN \rho_{\delta,1} \concANN \chi_1} \\
	&\leq \tfrac{81}{32} \kappa_1 \idconstANN^3 \CNdepth{\CNelt} \CNwidth{\CNelt}^{\kappa_2 + 6} \big[\!\log_2\!\big( \sqrt{\CNwidth{\CNelt}} \delta^{-1} \big)\big]^{\kappa_3}.
\end{split}
\end{equation*}
The rest of the proof of \cref{theorem_main} carries over without any changes.
\end{IEEEproof}


\section{Proofs for Section \ref{section_applications}}

\begin{IEEEproof}[\bf{Proof of \cref{ex_Lipschitz_sum_bounded}}]
Let $\Catalog = \Lipcatalog$ be the $\Lipconst$-Lipschitz catalog and suppose $\approxSet{} = (\approxSet{f})_{f \in \Catalog}$
and $\Lip = (\Lip_{f})_{f \in \Catalog}$ are defined as in \eqref{ExLip(i)} of \cref{ex_Lip_catalog}.
For $d \in \N$, let $V_d \in \R^{1 \times d}$ be the matrix $V_d = (1, \dots,1)$ and $\CNelt_d \in \CNarb$
the catalog network $\CNelt_d = [ (\ID{d},0,(f_1,\dots,f_d) ), (V_d,0,\ID{} ) ]$.
Then $\CNfct{\CNelt_d} = g_d$, $\CNdepth{\CNelt_d} = 2$, $\CNwidth{\CNelt_d} = d$ and $\CNtechLip{\CNelt_d} \leq \sqrt{d} \Lipconst$.
Moreover, $\CNredom{\CNelt_d} = \CNrdlayer{1}{\CNelt_d} \supseteq [-r,r]^d$ because $\CNrdlayer{2}{\CNelt_d} = \approxSet{\ID{}} = \R$.
Now, the proposition follows from \cref{rem_improved_exp_for_dim_one_thrm,cor_ReLU_weight_one}.
\end{IEEEproof}

\vspace{0.2cm}

\begin{IEEEproof}[\bf{Proof of \cref{ex_Lipschitz_sum_bounded_v2}}]
This proposition is proved as \cref{ex_Lipschitz_sum_bounded}, except that we fix $d$ in the beginning, use 
\eqref{ExLip(i)} of \cref{ex_Lip_catalog} with $d\Lipconst (r+1)$ instead of $r$ and define 
$\CNelt_d = [ (\ID{d},0,(f_1,\dots,f_d) ), (V_d,0,f_0 ) ]$.
Then, $\CNtechLip{\CNelt_d} \leq \sqrt{d} \Lipconst^2$ and $\CNrdlayer{2}{\CNelt_d} = [-d\Lipconst (r+1),d\Lipconst (r+1)]$, 
which ensures that $\CNredom{\CNelt_d} \supseteq [-r,r]^d$.
\end{IEEEproof}

\vspace{0.2cm}

\begin{IEEEproof}[\bf{Proof of \cref{ex_Lipschitz_max_bounded_v3}}]
	Let $\Catalog = \Lipmaxcatalog$ be the $\Lipconst$-Lipschitz-maximum catalog and suppose $\approxSet{} = (\approxSet{f})_{f \in \Catalog}$ and $\Lip = (\Lip_{f})_{f \in \Catalog}$ are defined as in \eqref{max(i)} of \cref{ex_Lip_max_catalog} 
(with $R = \Lipconst(r+1)$ instead of $r$).
Since we know that the functions $g_i$ are 1-Lipschitz, we may actually set $\Lip_{g_i} = 1$ for all $i \in \N$ without affecting the approximability of the catalog.
Consider $\CNelt_d \in \CNarb$, $d \in \N_{\geq 2}$, given by
\begin{equation*}
\begin{split}
	\CNelt_d = \big[ &\big(\ID{d},0,(f_1,\dots,f_d) \big), \\
	&\big(\ID{d},0,(\maxfct{2},\ID{},\dots,\ID{}) \big), \\
	&\big(\ID{d-1},0,(g_2,\ID{},\dots,\ID{}) \big), \\
	&\big(\ID{d-1},0,(\maxfct{2},\ID{},\dots,\ID{}) \big), \dots, \\
	&\big(\ID{2},0,(g_{d-1},\ID{}) \big), \big(\ID{2},0,\maxfct{2} \big), \\
	&\big(\ID{},0,g_d \big) \big].
\end{split}
\end{equation*}
Then $\CNfct{\CNelt_d} = h_d$, $\CNdepth{\CNelt_d} = 2d-1$, $\CNwidth{\CNelt_d} = d$ and $\CNtechLip{\CNelt_d} \leq \Lipconst$.
Moreover, $\CNrdlayer{2i}{\CNelt_d} = \R^{d-i+1}$ and $\CNrdlayer{2i+1}{\CNelt_d} = [-R,R] \times \R^{d-i-1}$ for all $i \in \{1,\dots,d-1\}$ as well as $\CNrdlayer{1}{\CNelt_d} = [-R,R]^d$.
If we define $H_i \colon \R^i \to \R$, $i \in \N_{\geq 2}$, by $H_i(x) = \max\{h_{i-1}(x_1,\dots,x_{i-1}),f_i(x_i)\}$, then it follows by induction that $H_i$ is $\Lipconst$-Lipschitz with respect to the sup-norm and $|H_i(0)| \leq \Lipconst$.
This proves $(H_i(x_1,\dots,x_i),x_{i+1},\dots,x_d) \in \CNrdlayer{2i-1}{\CNelt_d}$ for all $x \in [-r,r]^d$ and $i \in \{2,\dots,d\}$.
Since this corresponds to the evaluation of the first $2(i-1)$ layers of $\CNelt$, we have shown $[-r,r]^d \subseteq \CNredom{\CNelt_d}$.
Now we conclude with\footnote{Strictly speaking, \cref{rem_improved_exp_for_dim_one_thrm} is not applicable to the Lipschitz-maximum catalog, but we only used $\maxfct{2}$ and could remove $\maxfct{d}$, $d \geq 3$, from the catalog. This also enables us to use $\kappa = (\Lipconst,13\Lipconst^2r,0,1)$ instead of $\kappa_2 = 3$.} \cref{rem_improved_exp_for_dim_one_thrm,cor_ReLU_weight_one}.
\end{IEEEproof}

\vspace{0.2cm}

\begin{IEEEproof}[\bf{Proof of \cref{ex_Lipschitz_sum_unbounded}}]
Let $\Catalog = \Lipcatalog$ be the $\Lipconst$-Lipschitz catalog and suppose 
$\approxSet{} = (\approxSet{f})_{f \in \Catalog}$ and $\Lip = (\Lip_{f})_{f \in \Catalog}$ are defined as in 
\eqref{ExLip(ii)} of \cref{ex_Lip_catalog}. For all $d \in \N$, let 
$V_d \in \R^{1 \times d}$ be the matrix $V_d = (1, \dots, 1)$ and $\CNelt_d \in \CNarb$
the catalog network $\CNelt_d = [ (\ID{d},0,(f_1,\dots,f_d) ), (V_d,0,\ID{} ) ]$.
Then $\CNfct{\CNelt_d} = g_d$, $\CNtransl{\CNelt_d} = 1$, $\CNdepth{\CNelt_d} = 2$, $\CNwidth{\CNelt_d} = d$ 
and $\CNtechLip{\CNelt_d} \leq \sqrt{d} \Lipconst$. Moreover, $\CNredom{\CNelt_d} = \CNrdlayer{1}{\CNelt_d} = \R^d$ 
because $\CNrdlayer{2}{\CNelt_d} = \approxSet{\ID{}} = \R$. Now, the proposition follows from
\cref{rem_improved_exp_for_dim_one_thrm,cor_ReLU}.
\end{IEEEproof}

\vspace{0.2cm}

\begin{IEEEproof}[\bf{Proof of \cref{ex_Lipschitz_prod_v3}}]
Let $\Catalog = \Lipprodcatalog$ be the $\Lipconst$-Lipschitz-product catalog and suppose 
$\approxSet{} = (\approxSet{f})_{f \in \Catalog}$ and $\Lip = (\Lip_{f})_{f \in \Catalog}$ 
are defined as in \cref{ex_Lip_prod_catalog} (with $\nicefrac{1}{\sqrt{8}}$ instead of $r$ and $R=\nicefrac{1}{\sqrt{8}}$).
Since we know that the functions $g_i$ are 1-Lipschitz, we may actually set $\Lip_{g_i} = 1$ for all $i \in \N$ without affecting the approximability of the catalog.
Consider $\CNelt_d \in \CNarb$, $d \in \N_{\geq 2}$, given by
\begin{equation*}
\begin{split}
	\CNelt_d = \big[ &\big(\ID{d},0,(f_1,f_2,\dots,f_d) \big), \\
	&\big(\ID{d},0,(\prodfct,\ID{},\dots,\ID{}) \big), \\
	&\big(\ID{d-1},0,(g_2,\ID{},\dots,\ID{}) \big), \\
	&\big(\ID{d-1},0,(\prodfct,\ID{},\dots,\ID{}) \big), \dots, \\
	&\big(\ID{2},0,(g_{d-1},\ID{}) \big), \big(\ID{2},0,\prodfct \big), \\
	&\big(\ID{},0,g_d \big) \big].
\end{split}
\end{equation*}
Then $\CNfct{\CNelt_d} = h_d$, $\CNdepth{\CNelt_d} = 2d-1$, $\CNwidth{\CNelt_d} = d$ and $\CNtechLip{\CNelt_d} \leq \Lipconst$.
Moreover, $\CNrdlayer{2i}{\CNelt_d} = [-R,R]^2 \times \R^{d-i-1}$ and $\CNrdlayer{2i+1}{\CNelt_d} = [-R,R] \times \R^{d-i-1}$ for all $i \in \{1,\dots,d-1\}$ as well as $\CNrdlayer{1}{\CNelt_d} = [-R,R]^d$.
If we define $H_i \colon \R^i \to \R$, $i \in \N_{\geq 2}$, by $H_i(x) = h_{i-1}(x_1,\dots,x_{i-1}) f_i(x_i)$, then $H_i(0) = 0$ and it follows by induction that $|H_i(x)| \leq R^i$ and $|h_i(x)| \leq R^i$ for all $x \in [-r,r]^d$ and $i \geq 2$.
This shows that $[-r,r]^d \subseteq \CNredom{\CNelt_d}$, and we can conclude with 
\cref{rem_improved_exp_for_dim_one_thrm,cor_ReLU_weight_one}.
\end{IEEEproof}

\vspace{0.2cm}

\begin{IEEEproof}[\bf{Proof of \cref{ex_prod}}]
	Assume without loss of generality that $d \geq 3$, let $\Catalog = \Lipprod$ be the product catalog, and suppose 
	 $\approxSet{} = (\approxSet{f})_{f \in \Catalog}$ and $\Lip = (\Lip_{f})_{f \in \Catalog}$ are defined as in \cref{ex_prod_catalog}.
Moreover, let $\CNelt_d \in \CNarb$ be given by
\begin{equation*}
\begin{split}
	\CNelt_d = \big[ &\big(\ID{d},0,(\prodfct,\ID{},\dots,\ID{}) \big), \\
	&\big(\ID{d-1},0,(\prodfct,\ID{},\dots,\ID{}) \big), \dots, \\
	&\big(\ID{3},0,(\prodfct,\ID{}) \big), \big(\ID{2},0,\prodfct \big) \big].
\end{split}
\end{equation*}
Then $\CNfct{\CNelt_d} = f_d$, $\CNdepth{\CNelt_d} = d-1$, $\CNwidth{\CNelt_d} = d$ and $\CNtechLip{\CNelt_d} = 8^{\nicefrac{(d-1)}{2}} r^{d(d-1)}$.
Moreover, $\CNrdlayer{n}{\CNelt_d} = [-r^d,r^d]^2 \times \R^{d-n-1}$ for all $n \in \{1,\dots,d-1\}$.
Hence, the fact that for all $n \in \{1,\dots,d-1\}$ and $x \in [-r,r]^d$ we have $|\!\sprod{i=1}{n} x_i| \leq r^d$ ensures that $[-r,r]^d \subseteq \CNredom{\CNelt_d}$.
Now, we can conclude with \cref{rem_improved_exp_for_dim_one_thrm,cor_ReLU_weight_one_LOG} using the inequality
\begin{equation*}
	\log_2\!\big( \sqrt{d} (d-1) 8^{(d-1)/2} \epsilon^{-1} \big) \leq \tfrac{31}{18} d \log_2(d) \log_2(\epsilon^{-1})
\end{equation*}
if $r=1$ and the inequality
\begin{equation*}
\begin{split}
	\log_2\!\big( \sqrt{d} (d-1) 8^{(d-1)/2} &r^{d(d-1)} \epsilon^{-1} \big) \\
	&\leq \tfrac{67}{54} \log_2(r) d^2 \log_2(d) \log_2(\epsilon^{-1})
\end{split}
\end{equation*}
if $r \geq 2$, for which we note that $\log_2(d) \geq \frac{3}{2}$ since $d \geq 3$.
\end{IEEEproof}

\vspace{0.2cm}

\begin{IEEEproof}[\bf{Proof of \cref{ex_ridge}}]
	Let $\Catalog = \Lipcatalog$ be the $\Lipconst$-Lipschitz catalog and suppose 
	$\approxSet{} = (\approxSet{f})_{f \in \Catalog}$ and $\Lip = (\Lip_{f})_{f \in \Catalog}$ are defined as in 
	\eqref{ExLip(i)} of \cref{ex_Lip_catalog} (with $drS$ instead of $r$).
Let $\CNelt_d \in \CNarb$ be given by $\CNelt_d = (\theta_d^T,0,f )$, where $^T$ denotes transposition.
Then $\CNfct{\CNelt_d} = g_d$, $\CNdepth{\CNelt_d} = 1$, $\CNwidth{\CNelt_d} = d$ and $\CNtechLip{\CNelt_d} \leq \sqrt{d} \Lipconst S$.
Moreover, $[-r,r]^d \subseteq \{x \in \R^d \colon |\theta_d \cdot x| \leq drS \} \subseteq \CNredom{\CNelt_d}$ by 
the Cauchy--Schwarz inequality. So the proposition follows from 
\cref{rem_improved_exp_for_dim_one_thrm,cor_ReLU_weight_one}.
\end{IEEEproof}

\vspace{0.2cm}

\begin{IEEEproof}[\bf{Proof of \cref{ex_Gaussian_RBF}}]
	Let $\Catalog = \GaussianRBFcat$ be the Gaussian radial basis function catalog and suppose 
	 $\approxSet{} = (\approxSet{f})_{f \in \Catalog}$ and $\Lip = (\Lip_{f})_{f \in \Catalog}$ are defined as in 
	 \cref{ex_RBF_catalog} (with $r+S$ instead of $r$).
For all $d \in \N$,
let $U_d \in \R^{dN \times d}$ be the block-matrix $U_d = (\ID{d},\dots,\ID{d})^T$,
$b_d \in \R^{dN}$ the vector $(v_{d,1},\dots,v_{d,N})^T$,
$V_d \in \R^{N \times dN}$ the block-matrix with $\alpha_i(1,\dots,1) \in \R^{1 \times d}$ blocks on the $i$-th entry of 
the diagonal and 0 entries otherwise and $W \in \R^{1 \times N}$ the matrix $(u_1,\dots,u_N)$.
Moreover, let $\CNelt_d \in \CNarb$ be given by
\begin{equation*}
	\big[ \big( U_d,-b_d,(\sqfct,\dots,\sqfct) \big) , \big( V_d,0,(\expfct,\dots,\expfct) \big) , \big( W,0,\ID{} \big) \big].
\end{equation*}
Then $\CNfct{\CNelt_d} = f_d$, $\CNdepth{\CNelt_d} = 3$, $\CNwidth{\CNelt_d} = dN$ and $\CNtechLip{\CNelt_d} \leq 2 (r+S) \sqrt{d}NS^2$.
Moreover, $\CNrdlayer{1}{\CNelt_d} = [-(r+S),r+S]^{dN}$, $\CNrdlayer{2}{\CNelt_d} = [0,\infty)^N$ and $\CNrdlayer{3}{\CNelt_d} = \R$.
It follows that $[-r,r]^d \subseteq \CNredom{\CNelt_d}$, and the proposition follows from
\cref{rem_improved_exp_for_dim_one_thrm,cor_ReLU_weight_one}.
\end{IEEEproof}

\end{appendices}


\section*{Acknowledgment}
The authors are grateful to Philippe von Wurstemberger for helpful comments and suggestions. 
This work has been funded in part by the Swiss National Science Foundation Research Grant 175699
and by the Deutsche Forschungsgemeinschaft (DFG, German Research Foundation) under Germany's Excellence Strategy EXC 2044-390685587, 
Mathematics M{\"u}nster: Dynamics - Geometry - Structure.

\IEEEtriggeratref{33}

\bibliographystyle{IEEEtran}
\bibliography{}

\end{document}